\documentclass[12pt,english]{amsart}
\usepackage{amsthm,amssymb,latexsym,xcolor}
\usepackage{hyperref,enumerate,times,graphicx}

\newcommand{\Aset}{\mathbb{A}}
\newcommand{\Cset}{\mathbb{C}}
\newcommand{\Nset}{\mathbb{N}}

\newcommand{\Rset}{\mathbb{R}}
\newcommand{\Tset}{\mathbb{T}}
\newcommand{\Zset}{\mathbb{Z}}

\theoremstyle{plain}
\newtheorem{theorem}{Theorem}
\newtheorem{proposition}[theorem]{Proposition}
\newtheorem{lemma}[theorem]{Lemma}
\newtheorem{corollary}[theorem]{Corollary}
\newtheorem{definition}{Definition}
\theoremstyle{remark}
\newtheorem{nota}{Notation}
\newtheorem{remark}{Remark}
\newtheorem{example}{Example}

\newcommand{\rmi}{\mskip2mu{\rm i}\mskip1mu}

\newcommand{\Order}{\mathop{\rm O}\nolimits}
\newcommand{\rmd}{\,{\rm d}}
\newcommand{\rme}{{\rm e}}
\newcommand{\Graph}{\mathop{\rm graph}\nolimits}
\newcommand{\av}{{\boldsymbol{a}}}
\newcommand{\hv}{{\boldsymbol{h}}}
\newcommand{\xv}{{\boldsymbol{x}}}
\newcommand{\yv}{{\boldsymbol{y}}}
\newcommand{\alphav}{{\boldsymbol{\alpha}}}
\newcommand{\betav}{{\boldsymbol{\beta}}}
\newcommand{\psiv}{{\boldsymbol{\psi}}}

\newcommand{\llbracket}{[\![}
\newcommand{\rrbracket}{]\!]}
        
\begin{document}

\title{High-order persistence of resonant caustics in
perturbed circular billiards}

\author{Comlan Edmond Koudjinan}
\address{University of Toronto, Ontario, Canada}
\email{koudjinanedmond@gmail.com}
\author{Rafael Ram{\'\i}rez-Ros}
\address{Universitat Polit\`ecnica de Catalunya, Barcelona, Spain}
\email{rafael.ramirez@upc.edu}

\date{September 8, 2025
     (revised version,
      to appear in \emph{Ergodic Theory \& Dynamical Systems})}

\begin{abstract}
We find necessary and sufficient conditions for high-order
persistence of resonant caustics in perturbed circular billiards.
The main tool is a perturbation theory based
on the Bialy-Mironov generating function for convex billiards.
All resonant caustics with period $q$ persist
up to order $\lceil q/n \rceil -1$ under any
polynomial deformation of the circle of degree $n$.
\end{abstract}


\maketitle

\section{Introduction}

The goal of this work is two-fold.
First,
to extend the first-order perturbation theory for exact twist maps
developed in~\cite{RamirezRos2006,PintoRamirezRos2013,Damasceno_etal2017}
to a higher-order theory.
Second,
to apply that theory to the study of high-order persistence
of resonant caustics in perturbed circular billiards.
The second goal is strongly motivated by some of the numerical
experiments discussed in~\cite{MartinTamaritRamirezRos2016b}.

The computational aspects of our analysis are greatly simplified
when working with the Bialy-Mironov generating function for convex
billiards discovered in~\cite{BialyMironov2017,Bialy2018}.
We were also inspired by~\cite{BialyTabachnikov2022}.
We know just a few practical high-order Melnikov theories for
time-periodic perturbations of integrable continuous systems (ODEs)
---see, for instance, ~\cite{ChenWang2017,ChenWang2019}---,
but none for perturbations of integrable discrete systems (maps).
In that sense, our theory is novel.

The fragility of resonant caustics is a key idea behind recent
proofs of local versions of the Birkhoff conjecture (see below)
and related results about the rigidity of the length spectrum
of strictly convex
domains~\cite{AvilaDeSimoiKaloshin2016,KaloshinSorrentino2018b,
HuangKaloshinSorrentino2018,HuangKaloshin2019,KaloshinSorrentino2021,
HezariZelditch2022}.
See also the
surveys~\cite{KaloshinSorrentino2018a,FierobeKaloshinSorrentino2024}.
Almost all these works describe
the first-order persistence condition of resonant caustics
contained in~\cite{RamirezRos2006}.
We hope that our new high-order persistence conditions
will be equally useful.

A \emph{caustic} is a curve such that any billiard trajectory,
once tangent to the curve, stays tangent after every reflection.
The robustness of a convex caustic is closely related to the
arithmetic properties of its \emph{rotation number} $\rho \in (0,1)$,
a number that measures the number of turns around the caustic
per impact.
Tangent lines to the caustic can be counterclockwise or
clockwise oriented.
We fix the counterclockwise orientation, so $\rho \in (0,1/2]$.
Lazutkin~\cite{Lazutkin1973} showed that for
any smooth strictly convex domain there is a
positive measure Cantor set $\mathcal{R} \subset (0,1/2)$
of Diophantine rotation numbers that accumulates to $0$
such that there is a caustic for any rotation number
$\rho \in \mathcal{R}$.
These caustics persist under smooth deformations of the
domain~\cite{Popov1994}.

Let $\rho = p/q \in (0,1/2]$ be a rational rotation number
such that $\gcd(p,q) = 1$.
A convex caustic is called \emph{$p/q$-resonant}
(or $p/q$-\emph{rational}) when all its tangent trajectories
form closed polygons with $q$ sides that make $p$ turns around the caustic.
We say that $q$ is the \emph{period} of the caustic.
Resonant caustics generically break up under perturbation.
Recent results in~\cite{FierobeSorrention2024} confirm their
fragility.
Once fixed $q \ge 2$,
the space of convex domains with a resonant caustic of period $q$ has
infinite dimension and codimension~\cite{BaryshnikovZharnitsky2006}.
The space of convex domains with at least one resonant
caustic is dense in the space of all convex domains~\cite{KaloshinZhang2018}.

We shall not deal with the case $\rho = 1/2$,
since convex domains with $1/2$-resonant caustics
are easily characterized as the \emph{constant width}
domains~\cite{Knill1998,Gutkin2012}.
Centrally symmetric convex domains with a $1/4$-resonant
caustic have also been completely characterized in terms of
the Fourier coefficients of the square of the support function
of the convex domain in~\cite{BialyTsodikovich2023a}.
Some non-circular convex domains with a $1/3$-resonant
caustic were constructed in~\cite{Innami1988}.

Circles and ellipses are the only known
strictly convex smooth domains almost completely foliated
by convex caustics.
The centenary \emph{Birkhoff conjecture} claims that they are
the only ones~\cite{Tabachnikov1995}.
Bialy~\cite{Bialy1993} proved the following weak version of this conjecture.
If almost every billiard trajectory in a convex domain
is tangent to a convex caustic, then the domain is a disk.
A much stronger version of the Birkhoff conjecture for
centrally symmetric $C^2$-domains,
based on the structure of the $1/4$-resonant caustic,
was recently established by Bialy and Mironov~\cite{BialyMironov2022}.
See also~\cite{BialyTsodikovich2023a,Bialy2024} for effective
(that is, quantitative) versions on these two results.
Near centrally symmetric domains were considered
in~\cite{KaloshinKoudjinanZhang2024}.

We are interested in two practical problems.
First, to characterize the deformations of the circle
that preserve a given resonant caustic.
Second, to determine all resonant caustics that are preserved
under a given deformation of the circle.
In that regard,
we recall that any $\Zset_2$-symmetric analytical deformation of
a circle (with certain Fourier decaying rate)
preserving both its $1/2$-resonant and $1/3$-resonant caustics
has to be an isometric transformation~\cite{Zhang2019}.

In what follows we introduce some notations and
state our two main results.

Let $\Gamma_\epsilon$ be a deformation of the unit circle
with smooth support function
\begin{equation}
\label{eq:SupportFunction}
h(\psi;\epsilon) =
h_\epsilon(\psi) \asymp
1 + \sum_{k \ge 1} \epsilon^k h_k(\psi)
\quad \mbox{ as } \epsilon \to 0,
\end{equation}
where $\psi \in \Tset = \Rset/2\pi \Zset$ is the normal angle and
$\epsilon \in [-\epsilon_0,\epsilon_0]$ is the perturbative parameter.
We say that a resonant caustic of the unit circle
\emph{$\Order(\epsilon^m)$-persists} under $\Gamma_\epsilon$
when the billiard in $\Gamma_\epsilon$ is
$\Order(\epsilon^{m+1})$-close to having that resonant caustic.
See Definition~\ref{defi:PersistenceOrder} for more details.

Let $\nu_l:(0,1/2) \to \Rset \cup \{\infty\}$,
with $l \in \Zset$ and $|l| \ge 2$,
be the sequence of functions given by
\begin{equation}
\label{eq:Factors}
\nu_l(\rho) = \nu_{-l}(\rho) =
\begin{cases}
\displaystyle
\frac{\tan(l \pi \rho) - l \tan (\pi \rho)}{\tan (\pi \rho) \tan(l\pi \rho)}, &
\mbox{if $2 l \rho \not \in \Zset$}, \\
1/\tan(\pi \rho), &
\mbox{if $2l \rho \in \Zset$ but $l \rho \not \in \Zset$}, \\
\infty, & \mbox{if $l \rho \in \Zset$}.
\end{cases}
\end{equation}

Once again, we realize that Gutkin's equation
$\tan(l\pi \rho) = l \tan(\pi \rho)$
is ubiquitous in billiard problems.
See~\cite{Gutkin2012,Bialy2018,Bialy_etal2019,BialyTsodikovich2023b}
for other examples.
Cyr~\cite{Cyr2012} proved that $\nu_l(\rho)$
has no rational roots $\rho = p/q \in (0,1/2)$ when $|l| \ge 2$.
The case $l \rho \in \Zset$ never takes place in our computations.
The singular value $\nu_l(\rho) = \infty$ has been written
just for definitness.
It is irrelevant.

Fourier coefficients of $2\pi$-periodic functions are denoted
with a hat: $a(t) = \sum_{l \in \Zset} \hat{a}_l \rme^{\rmi l t}$.
Given a $2\pi$-periodic smooth function $a(t)$ and
a subset $R \subset \Zset$, let
$\mu_R\{a(t)\} =
 \sum_{l \in R} \hat{a}_l \rme^{\rmi l t}$
be the projection of $a(t)$ onto its $R$-harmonics.
We only consider the cases
$R = q\Zset$ or $R = q\Zset^\star$
with $\Zset^\star = \Zset \setminus \{ 0 \}$.

\begin{theorem}
\label{thm:Persistence} 
Let $\rho = p/q \in (0,1/2)$ be any rational rotation number
such that $\gcd(p,q) = 1$.
The high-order persistence of the $p/q$-resonant caustic of
the unit circle under the deformation with support
function~\eqref{eq:SupportFunction} can be determined
as follows.
\begin{enumerate}[a)]
\item
It $\Order(\epsilon)$-persists if and only if
$\mu_{q\Zset^\star} \{ h_1 \} = 0$.

\item
It $\Order(\epsilon^2)$-persists if and only if it
$\Order(\epsilon)$-persists and 
$\mu_{q\Zset^\star} \{ h_2 + \theta^2_1/2 \} = 0$,
where
\[
\theta_1(t) =
\sum_{l \not \in q\Zset \cup \{-1,1\}}
\nu_l(\rho) \hat{h}_{1,l} \ \rme^{\rmi l t}
\quad \mbox{ if } \quad
h_1(\psi) =
\sum_{l \in \Zset} \hat{h}_{1,l} \ \rme^{\rmi l \psi}.
\]

\item
It $\Order(\epsilon^m)$-persists for some $m \ge 3$ if and only if
it $\Order(\epsilon^{m-1})$-persists and
$\mu_{q\Zset^\star} \{ h_m + \zeta_m \} = 0$,
where $\zeta_m$ is a smooth $2\pi$-periodic function,
only depending on $h_1,\ldots,h_{m-1}$,
that can be explicitly computed from recurrences given
along the paper.
\end{enumerate}
\end{theorem}

The $\Order(\epsilon)$-persistence result in Theorem~\ref{thm:Persistence}
is just a reformulation of the main theorem in~\cite{RamirezRos2006}.
Condition $\mu_{q\Zset^\star} \{ h_m + \zeta_m \} = 0$
is equivalent to $\mu_{q\Zset}\{ h'_m + \zeta'_m\} = 0$.
In particular,
condition $\mu_{q\Zset^\star} \{ h_2 + \theta_1^2/2 \}$
is equivalent to $\mu_{q\Zset} \{ h'_2 + \theta_1 \theta'_1 \} = 0$.

Let $T_n[\psi]$ be the space of $2\pi$-periodic real
trigonometric polynomials of degree $\le n$ in $\psi$.

\begin{definition}
\label{defi:PolynomialDeformation}
We say that a deformation $\Gamma_\epsilon$ of the unit circle
with support function~\eqref{eq:SupportFunction} is
\emph{polynomial} of \emph{degree} $\le n$ when
\begin{equation}
\label{eq:DegreeCondition}
h_k(\psi) \in T_{n k}[\psi], \quad \forall k \ge 1,
\end{equation}
and is \emph{centrally} or
\emph{anti-centrally symmetric} when
$h_\epsilon(\psi+\pi) = h_\epsilon(\psi)$
or $h_\epsilon(\psi+\pi) = h_{-\epsilon}(\psi)$.
\end{definition}

Being centrally symmetric is a property of single curves.
Being anti-centrally symmetric is a property of deformations.

\begin{theorem}
\label{thm:Polynomial}
Let $\lceil x \rceil = \min\{ k \in \Zset: k \ge x \}$
be the ceil function.
If $\rho = p/q \in (0,1/2)$ is a rational rotation number such that
$\gcd(p,q) = 1$ and $\Gamma_\epsilon$ is a polynomial deformation
of the unit circle of degree $\le n$,
then the $p/q$-resonant caustic $\Order(\epsilon^{\chi -1})$-persists
under $\Gamma_\epsilon$, where
\[
\chi = \chi(\Gamma_\epsilon,q) =
\begin{cases}
1 + 2 \left\lceil (q-n)/2n \right\rceil, &
\text{for anti-centrally symmetric $\Gamma_\epsilon$ and odd $q$,} \\
2\lceil  q/2n \rceil, &
\mbox{for anti-centrally symmetric $\Gamma_\epsilon$ and even $q$}, \\
 \lceil 2q/n  \rceil, &
 \mbox{for centrally symmetric $\Gamma_\epsilon$ and odd $q$}, \\
 \lceil  q/n  \rceil, & \mbox{otherwise.}
\end{cases}
\]
\end{theorem}

The idea behind this theorem is quite simple.
For non-symmetric deformations,
it suffices to check that $\zeta_m(t) \in T_{nm}[t]$ for
$m = 1,\ldots,\chi-1$, where $\zeta_m(t)$ are the functions
introduced in Theorem~\ref{thm:Persistence}.
Symmetric deformations require to check that, in addition,
those polynomials $\zeta_m(t)$ are
$\pi$-periodic or $\pi$-antiperiodic.

Polynomial deformations of the unit circle of degree $\le n$
can be defined without mentioning support functions.
For instance, we can define them in Cartesian coordinates $(x,y)$ as
\begin{equation}
\label{eq:PolynomialPerturbation}
\Gamma_\epsilon =
\big\{
(x,y) \in \Rset^2 : x^2 + y^2 = P(x,y;\epsilon)
\big\},
\end{equation}
for some smooth function $P(x,y;\epsilon)$ of the form
$P(x,y;\epsilon) \asymp 1 + \sum_{k \ge 1} \epsilon^k P_k(x,y)$
as $\epsilon \to 0$ with
$P_k(x,y) \in \Rset_{kn}[x,y]$ for all $k \ge 1$.
Alternatively,
we can also define them in polar coordinates $(r,\phi)$ as
\[
\Gamma_\epsilon =
\big\{
r(\phi;\epsilon) \cdot (\cos \phi,\sin \phi) :
\phi \in \Tset
\big\},
\]
for some smooth \emph{polar function} $r(\phi;\epsilon)$
of the form
$r(\phi;\epsilon) \asymp 1 + \sum_{k\ge 1} \epsilon^k r_k(\phi)$
as $\epsilon \to 0$ with
$r_k(\phi) \in T_{kn}[\phi]$ for all $k \ge 1$.
The Cartesian setting was considered 
in~\cite{MartinTamaritRamirezRos2016b} with
$P(x,y;\epsilon) = 1 - \epsilon y^n$.
The polar setting was considered in~\cite{RamirezRos2006}
with $r(\phi;\epsilon) = 1 +\epsilon r_1(\phi) + \Order(\epsilon^2)$
and in~\cite{Zhang2019} with
$r(\phi;\epsilon) =
 1 +\epsilon r_1(\phi) + \epsilon^2 r_2(\phi) + \Order(\epsilon^3)$.
Any deformation of the unit circle expressed in
Cartesian coordinates as~\eqref{eq:PolynomialPerturbation}
for some $P(x,y;\epsilon) = 1 + \epsilon P_1(x,y)$ with
$P_1(x,y) \in \Rset_n[x,y]$ is a polynomial deformation
of degree $\le n$ in the sense of
Definition~\ref{defi:PolynomialDeformation}
with $h_1(\psi) = \frac{1}{2} P_1(\cos \psi,\sin \psi)$.
See Lemma~\ref{lem:Perturbation}.
We are interested in such deformations because we want
to understand the numerical experiments
discussed in~\cite{MartinTamaritRamirezRos2016b}.
However, for brevity, we omit the corresponding proofs for
deformations written in polar coordinates or in Cartesian
coordinates with more than the first order term $\epsilon P_1(x,y)$.
Such proofs are just a slew of boring computations based in the
Taylor, multinomial and Lagrange inversion theorems.
We only stress that $h_1(\psi) = r_1(\psi)$,
which justifies that the $\Order(\epsilon)$-persistence result
in Theorem~\ref{thm:Persistence} is just a reformulation of
the main theorem in~\cite{RamirezRos2006}.
 
The map $q \mapsto \chi(\Gamma_\epsilon,q)$ is unbounded for any
polynomial deformation $\Gamma_\epsilon$ of degree $\le n$, since
$\chi(\Gamma_\epsilon,q) \asymp 2q/n$ as odd $q \to +\infty$
for centrally symmetric deformations,
and $\chi(\Gamma_\epsilon,q)  \asymp q/n$ as $q \to +\infty$
otherwise.
The experiments described
in~\cite[Numerical Result 5]{MartinTamaritRamirezRos2016b},
which cover degrees $3\le n\le 8$ and periods $3\le q \le 100$,
suggest that none of the $p/q$-resonant caustics
$\Order(\epsilon^\chi)$-persists
under monomial deformations~\eqref{eq:PolynomialPerturbation}
with $P(x,y;\epsilon) = 1 - \epsilon y^n$ and $n \ge 3$.
Its proof requires to check that
$\mu_{q\Zset^\ast}\{h_\chi + \zeta_\chi \} \neq 0$,
which is a challenge.
If $\Gamma_\epsilon$ is a polynomial deformation
of the unit circle of degree $\le n$ such that
$\mu_{q\Zset^\ast}\{h_\chi + \zeta_\chi \} \neq 0$
for all $p/q \in (0,1/2)$,
then $\Gamma_\epsilon$ \emph{would break all resonant caustics}
in such a way that \emph{there would be breakups of any order},
because the map $q \mapsto \chi(\Gamma_\epsilon,q) \in \Nset$
is exhaustive.

The paper is organized as follows.
Section~\ref{sec:Existence} begins with a description of the
Bialy-Mironov generating function and ends with a list of
necessary conditions for the existence of smooth convex resonant
caustics in smooth strictly convex domains.
The general notion of high-order persistence of
convex resonant caustics in deformed smooth convex domains
is presented in Section~\ref{sec:Persistence}
and applied to deformed circles in Section~\ref{sec:Circles},
where Theorem~\ref{thm:Persistence} is proved.
The results about polynomial deformations of circles,
including Theorem~\ref{thm:Polynomial},
are presented in Section~\ref{sec:Polynomial}.
Finally, we discuss three open problems:
the co-preservation of resonant caustics with different
rotation numbers,
the convergence of a procedure to correct the original
deformation in order to preserve a chosen resonant caustic
and the asymptotic measure of some exponentially small phenomena
as the period $q$ grows.
See Section~\ref{sec:OpenProblems}.
Several technical proofs have been relegated to the appendices.

\section{Existence of smooth convex resonant caustics}
\label{sec:Existence}

To begin with,
we introduce coordinates in the space of oriented lines,
define the support function and the billiard map of a convex domain,
and describe the Bialy-Mironov generating function
following~\cite{BialyTabachnikov2022}.
We also recall the periodic version of the variational
principle for twist maps following~\cite{Meiss1992}.
Next, we combine all those elements to find necessary
conditions for the existence of smooth convex resonant caustics
in Theorem~\ref{thm:ExistenceConditions}.
This part is inspired by the
computations in~\cite[Theorem 2.2]{BialyTabachnikov2022}
and the Lagrangian approach to the existence of
rotational invariant curves (RICs) of twist maps described
in~\cite{LeviMoser2001,KaloshinZhang2018}.
Finally, we discuss five simple examples:
circles, ellipses, constant width curves, Gutkin billiard tables
(also called constant angle curves)
and centrally symmetric curves with $1/4$-resonant caustics.

The billiard dynamics acts on the subset of oriented lines (rays)
that intersect the boundary of the convex domain.
An oriented line $\ell$ can be written as
\[
\cos \varphi \cdot x + \sin \varphi \cdot y = \lambda,
\]
where $\varphi \in \Tset := \Rset/2\pi\Zset$ is the direction of
the right normal to the oriented line and $\lambda \in \Rset$
is the signed distance to the origin.
Thus, $(\varphi,\lambda) \in \Aset := \Tset \times \Rset$
are coordinates in the space of oriented lines,
which is topologically a cylinder.

Let $\Gamma$ be a smooth strictly convex closed curve of $\Rset^2$.
We fix its counterclockwise orientation and
assume its interior contains the origin $O$,
so there is a positive smooth $2\pi$-periodic function $h(\varphi)$,
called the \emph{support function} of $\Gamma$, such that
$\{ \lambda = h(\varphi) \}$ and $\{ \lambda = -h(\varphi+\pi) \}$
are the 1-parameter families of oriented lines
positively and negatively tangent to $\Gamma$.
Then
\[
z: \Tset \to \Gamma \subset \Rset^2 \simeq \Cset, \qquad
z(\varphi) =
\big( x(\varphi), y(\varphi) \big) =
h(\varphi) \rme^{\rmi \varphi} +
h'(\varphi) \rme^{\rmi (\varphi+\pi/2)},
\]
is a parametrization of $\Gamma$, where
$\varphi \in \Tset$ is the counterclockwise angle between
the positive $x$-axis and the outer normal to $\Gamma$ at
the point $z(\varphi)$.

The space of the oriented lines that intersect the interior
of $\Gamma$ is the open cylinder
\[
\Aset_\Gamma =
\left\{
(\varphi,\lambda) \in \Aset : -h(\varphi+\pi) < \lambda < h(\varphi)
\right\}
\]
and the \emph{billiard map} $f: \Aset_\Gamma \to \Aset_\Gamma$ acts
by the reflection law in $\Gamma$.
That is,
$f(\varphi,\lambda) = (\varphi_1,\lambda_1)$
means that the oriented line $\ell_1$ with coordinates
$(\varphi_1,\lambda_1)$ is the reflection of the oriented line $\ell$
with coordinates $(\varphi,\lambda)$ with respect to the tangent
to $\Gamma$ at the second intersection of $\ell$ with $\Gamma$.
See Figure~1.
The shocking discovery by Bialy and Mironov was that
the billiard map $f$ is an exact twist map with generating function
\[
S(\varphi,\varphi_1) = 2 h(\psi) \sin \theta, \quad
\psi = \frac{\varphi_1 + \varphi}{2}, \quad
\theta = \frac{\varphi_1 - \varphi}{2}.
\]
To be precise,
$\lambda_1 \rmd \varphi_1 - \lambda \rmd \varphi =
 f^* (\lambda \rmd \varphi) - \lambda \rmd \varphi = \rmd S$, so
\begin{equation}
\label{eq:ImplicitEquations}
f(\varphi,\lambda) = (\varphi_1,\lambda_1) \Leftrightarrow
\left\{
\begin{array}{ll}
\lambda  = -\partial_1 S(\varphi,\varphi_1) =
h(\psi) \cos \theta - h'(\psi) \sin \theta, \\
\lambda_1 = \partial_2 S(\varphi,\varphi_1) =
h(\psi) \cos \theta + h'(\psi) \sin \theta,
\end{array}
\right.
\end{equation}
and $f$ preserves the standard area form
$\rmd \varphi \wedge \rmd \lambda$.
See, for instance,~\cite[Proposition 2.1]{BialyTabachnikov2022}.
Here, $\partial_i S$ denotes the derivative with respect to the
$i$-th variable.
The strict convexity of $\Gamma$ implies the twist condition:
$\partial_{12}S(\varphi,\varphi_1) =
 \frac{1}{2} \rho(\psi)\sin \theta > 0$,
where $\rho(\psi) = h''(\psi) + h( \psi)$ is the
\emph{radius of curvature} of $\Gamma$ at the point $z(\psi)$.

\begin{figure}[t]
\begin{center}
\includegraphics*[width=10cm]{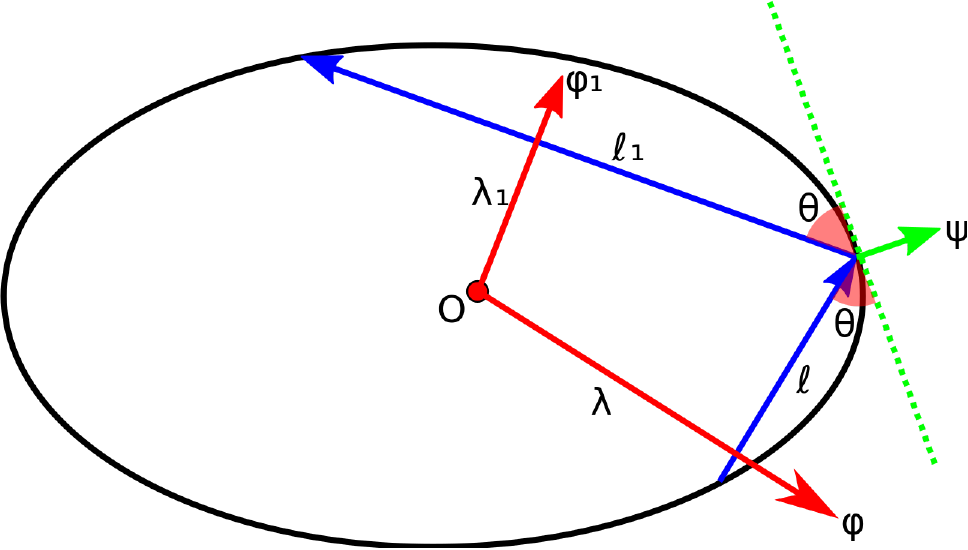}
\end{center}
\caption{\small
The billiard map $f(\varphi,\lambda) = (\varphi_1,\lambda_1)$,
the normal angle $\psi = (\varphi_1 + \varphi)/2$ and the
incidence-reflection angle $\theta = (\varphi_1 - \varphi)/2$.}
\label{fig:BilliardMap}
\end{figure}

We say that
$\psi \in \Tset$ and $\theta \in (0,\pi)$ are
the \emph{normal angle} and the \emph{angle of incidence/reflection}
at each impact point,
whereas $\varphi \in \Tset$ is the \emph{side angle}.
We consider $S(\varphi,\varphi_1)$ defined on
the universal cover
$\{ (\varphi,\varphi_1) \in \Rset^2 :
    \varphi < \varphi_1 < \varphi + 2\pi \}$
because $S(\varphi+2\pi,\varphi_1+2\pi) = S(\varphi,\varphi_1)$.

Let us briefly recall the classical variational principle
for exact twists maps.
The interested reader can find more details
in~\cite[\S V]{Meiss1992}.
We will introduce several operators which
are not standard in the variational approach,
but they simplify computations and shorten formulas.

Given any sequence $\{a_j\}$,
we consider the \emph{shift}, \emph{sum}, \emph{difference}
and \emph{$q$-average} operators
\[
\tau\{ a_j \} = a_{j+1},\quad
\sigma \{ a_j \} = a_{j+1} + a_j,\quad
\delta \{ a_j \} = a_{j+1} - a_j, \quad
\mu \{ a_j \} = \frac{a_j + \cdots + a_{j+q-1}}{q}.
\]
For simplicity, we omit the dependence of $\mu$ on $q$ and
we just say that $\mu$ is the \emph{average} operator.

Let $p/q \in (0,1)$ be a rational number with $\gcd(p,q) = 1$.
We call the elements of the space
\[
X =
\left\{
\{\varphi_j \} \in \Rset^\Zset :
\varphi_j < \varphi_{j+1} < \varphi_j + 2\pi, \
\varphi_{j+q} = \varphi_j + 2\pi p, \
\forall j \in \Zset
\right\}
\]
\emph{$p/q$-periodic sequences}.
We define the \emph{$p/q$-periodic action} $A: X \to \Rset$ as
\[
A\{ \varphi_j \} =
\sum_{j=0}^{q-1} S(\varphi_j,\varphi_{j+1}) =
2 \sum_{j=0}^{q-1} h(\psi_j) \sin \theta_j,
\]
where $\psi_j = \sigma \{ \varphi_j \}/2$
and $\theta_j = \delta\{ \varphi_j \}/2$.
Periodicities
$\varphi_{j+q} = \varphi_j + 2\pi p$,
$\psi_{j+q} = \psi_j + 2\pi p$ and $\theta_{j+q} = \theta_j$
imply that
$\mu \{ S(\varphi_{j_0+j},\varphi_{j_0+j+1}) \} =
 2\mu\{ h(\psi_{j_0+j}) \sin \theta_{j_0+j} \} =
 A\{ \varphi_j \}/q$ for all $j_0 \in \Zset$.
Critical points of the action,
which we call \emph{$p/q$-periodic configurations},
can be lifted to full \emph{$p/q$-periodic orbits}
of the billiard map $f$ by taking
$\lambda_j = -\partial_1 S(\varphi_j,\varphi_{j+1}) =
 \partial_2 S(\varphi_{j-1},\varphi_j)$.
Thus,
any $p/q$-periodic configuration defines a
$p/q$-periodic billiard trajectory inside $\Gamma$
with side angles $\varphi_j$,
normal angles $\psi_j = \sigma \{ \varphi_j \}/2$
and incidence-reflection angles $\theta_j = \delta\{ \varphi_j \}/2$.

We characterize such configurations in the next proposition,
where we also recall a formula for the length
of their corresponding periodic billiard trajectories.

\begin{proposition}
\label{prop:PeriodicConfiguration}
A $p/q$-periodic sequence $\{ \varphi_j \}$
is a $p/q$-periodic configuration if and only if
\[
\sigma \{ h'(\psi_j) \sin \theta_j \} -
\delta \{ h(\psi_j) \cos \theta_j \} = 0,
\quad \forall j \in \Zset,
\]
where $\psi_j = \sigma \{ \varphi_j \}/2$
and $\theta_j = \delta\{ \varphi_j \}/2$,
in which case
\[
2q \mu\big \{ h(\psi_j) \sin \theta_j \big \} =
A\{ \varphi_j\} = L,\quad
\mu \{ h'(\psi_j) \sin \theta_j \} = 0, \quad
\mu \{\theta_j\} = \pi p/q,
\]
where $L$ is the length of the $p/q$-periodic configuration
$\{ \varphi_j \}$.
\end{proposition}

\begin{proof}
A $p/q$-periodic sequence $\{ \varphi_j \}$ is a critical point
of the action if and only if
\[
\sigma \{ h'(\psi_j) \sin \theta_j \} -
\delta \{ h(\psi_j) \cos \theta_j \} =
\partial_2 S(\varphi_j,\varphi_{j+1}) +
\partial_1 S(\varphi_{j+1},\varphi_{j+2}) = 0
\]
for all $j \in \Zset$.
The first equality above follows from
$S(\varphi_j,\varphi_{j+1}) = 2 h(\psi_j) \sin \theta_j$
and relations
$\psi_j =
 \sigma \{ \varphi_j \}/2 =
 (\varphi_{j+1} + \varphi_j)/2$ and
$\theta_j =
 \delta \{ \varphi_j \}/2 =
 (\varphi_{j+1} - \varphi_j)/2$.

Identities
$2q \mu\big \{ h(\psi_j) \sin \theta_j \big \} = L$ and
$\mu \{ h'(\psi_j) \sin \theta_j \} = 0$ are proved
in~\cite[Theorem~2.2]{BialyTabachnikov2022}.
By $q$-periodicity of $\{\theta_j\}$, we get
\[
\mu \{ \theta_j \} =
\frac{1}{q} \sum_{i=0}^{q-1} \theta_{j+i} =
\frac{1}{q} \sum_{i=0}^{q-1} \frac{\theta_{j+i} + \theta_{j+i+1}}{2} =
\frac{1}{q} \sum_{i=0}^{q-1} \frac{\psi_{j+i+1} - \psi_{j+i}}{2} =
\frac{\psi_{j+q} - \psi_j}{2q} = \frac{\pi p}{q}.
\qedhere
\]
\end{proof}

We look for necessary conditions for the existence of
convex resonant caustics inside $\Gamma$.
Tangent lines to such a caustic can be counterclockwise or
clockwise oriented.
We fix the counterclockwise orientation,
so we assume that $p/q < 1/2$ from now on.
Any $p/q$-resonant convex caustic gives rise to a
1-parameter family of $p/q$-periodic orbits of the billiard map.
We want to parametrize this family using a
\emph{dynamical parameter} $t \in \Rset$ in which
the billiard map acts as the constant shift
$t \mapsto t + \omega$ with angular frequency $\omega = 2\pi p/q$.
The dynamical parameter is not unique.
If $a(t)$ is any smooth $\omega$-periodic function such that
$1 + a'(t) > 0$,
then $s = t + a(t)$ is another dynamical parameter.

Let us stress the three main differences between this setting,
where we deal with functions of a continuous variable $t \in \Rset$,
and the previous setting, where we had sequences whose elements
are labeled by a discrete index $k \in \Zset$.

Firstly,
we define the \emph{shift}, \emph{sum},
\emph{difference} and \emph{average} operators as
\[
\tau\{ a(t) \} = a(t + \omega),\quad
\sigma\{ a(t) \} = a(t + \omega) + a(t), \quad
\delta\{ a(t) \} = a(t + \omega) - a(t)
\]
and $\mu \{ a(t) \} = \frac{1}{q} \sum_{j=0}^{q-1} a(t + j \omega)$.
These operators diagonalize in the Fourier basis.
Operator $\mu$ is the projection onto the resonant
$q\Zset$-harmonics: $\mu = \mu_{q\Zset}$,
but we omit the $q\Zset$ subscript for simplicity.
Both claims are proved in Appendix~\ref{app:Operators}.

Secondly, we define the \emph{$p/q$-periodic action}
of a \emph{side function} $\varphi(t)$ as
\[
A\{ \varphi(t) \} =
q \mu \left\{ S(\varphi(t),\varphi(t+\omega) ) \right\} =
2q \mu \{ h(\psi(t)) \cdot \sin \theta(t) \},
\]
where $\varphi(t)$,
the \emph{normal function} $\psi(t)$ and
the \emph{incidence-reflection function} $\theta(t)$ are related by
\begin{equation}
\label{eq:Relations}
\psi =  \varphi + \theta = \sigma\{ \varphi \}/2, \quad
2\theta = \delta \{ \varphi \}, \quad
\mu\{ \psi - t \} = 0, \quad
\mu\left\{ t - \varphi \right\} = \omega/2 =\mu\{ \theta \}.
\end{equation}
We ask functions $\varphi(t) - t$, $\psi(t) - t$ and
$\theta(t)$ to be $2\pi$-periodic,
as a continuous analogue of the discrete periodicity conditions
$\varphi_{j+q} = \varphi_j + 2\pi p$,
$\psi_{j+q} = \psi_j + 2\pi p$ and $\theta_{j+q} = \theta_j$.
Hence, they are lifts of some functions
$\varphi,\psi: \Tset \to \Tset$ and
$\theta: \Tset \to \Rset$.
To simplify the exposition, we sometimes abuse the notation
and use the same symbol for an object and its lift.
Other times we denote lifts with a tilde.
We also ask that $\varphi'(t) > 0$,
so that $\varphi:\Rset \to \Rset$ can be inverted.

Thirdly, condition
$\sigma \{ h'(\psi_j) \sin \theta_j \} -
 \delta \{ h(\psi_j) \cos \theta_j \} = 0$
becomes the difference equation
\begin{equation}
\label{eq:DynamicEquation}
\sigma \{ h'\circ \psi \cdot \sin \theta \} -
\delta \{ h \circ \psi \cdot \cos \theta \} = 0.
\end{equation}

\begin{remark}
\label{rem:Relations}
Relations~\eqref{eq:Relations} are redundant,
but we have listed all of them for future references.
They can be used to determine
all three functions $\varphi(t)$, $\psi(t)$ and $\theta(t)$
from any one of them.
Usually, we will determine $\varphi(t)$ and $\psi(t)$ from $\theta(t)$.
If $\theta(t)$ is a smooth $2\pi$-periodic function
such that $\mu\{ \theta \} = \omega/2$,
then there is a unique smooth $2\pi$-periodic function
$\varphi(t)-t$ such that $\delta\{ \varphi \} = 2\theta$
and $\mu\{ t - \varphi \} = \omega/2$.
See Lemma~\ref{lem:Inversion} in Appendix~\ref{app:Operators}
for a proof.
Then $\psi = \varphi + \theta$ implies that
$\psi = \sigma\{ \varphi \}/2$ and $\mu\{ \psi - t \} = 0$.
\end{remark}

\begin{theorem}
\label{thm:ExistenceConditions}
Let $p/q \in (0,1/2)$ be any rational rotational number 
with $\gcd(p,q) = 1$.
If there is a smooth convex $p/q$-resonant caustic inside $\Gamma$,
the following necessary conditions hold.
\begin{enumerate}
\item[{\rm ({\bf f})}]
There are three smooth $2\pi$-periodic {\bf f}unctions
$\varphi(t) - t$, $\psi(t) - t$ and $\theta(t)$
related by~\eqref{eq:Relations} that satisfy the difference
equation ~\eqref{eq:DynamicEquation}.
Besides, $\varphi'(t) > 0$.

\item[{\rm ({\bf p})}]
There is a smooth {\bf p}arametrization $c:\Tset \to \Aset_\Gamma$
such that $G = c(\Tset)$ is a graph and $f(c(t)) = c(t+\omega)$.

\item[{\rm ({\bf a})}]
The $p/q$-periodic {\bf a}ction is constant on the side function:
$A\{\varphi(t)\} = L$, where $L$ is the length
of all $p/q$-periodic billiard trajectories in $\Gamma$.
\end{enumerate}
\end{theorem}

\begin{proof}
\begin{enumerate}
\item[{\rm ({\bf f})}]
Let $g(\varphi)$ be the smooth support function of the
$p/q$-resonant caustic.
Set
\[
G = \Graph(g)
  := \{ (\varphi,\lambda) \in \Aset : \lambda = g (\varphi) \big\}.
\]
The caustic is inside $\Gamma$, so $G \subset \Aset_\Gamma$.
Clearly, $G$ is $f$-invariant,
so $f|_G$ defines a smooth preserving orientation
circle diffeomorphism $r: \Tset \to \Tset$ such that
\[
\tilde{f}\big( \varphi,\tilde{g}(\varphi) \big) =
\big( \tilde{r}(\varphi),\tilde{g}(\tilde{r}(\varphi)) \big), \qquad
\tilde{r}^q(\varphi) = \varphi + 2\pi p = \varphi + q\omega,
\]
where $\tilde{f}:\tilde{\Aset}_\Gamma \to \tilde{\Aset}_\Gamma$,
$\tilde{g}:\Rset \to \Rset$ and $\tilde{r}:\Rset \to \Rset$ are
the corresponding lifts.
Identity $\tilde{r}^q(\varphi) = \varphi + q\omega$ is the key.
It follows from the definition of $p/q$-resonant caustic.
It implies that $\varphi \mapsto \tilde{r}(\varphi)$
becomes the constant shift $t \mapsto t +\omega$
in the smooth parameter
\[
t = \tilde{s}(\varphi) :=
\frac{1}{q}
\sum_{j=0}^{q-1} \big( \tilde{r}^j(\varphi) - j \omega \big).
\]
Inversion of $\tilde{s}(\varphi(t)) = t$ defines a
smooth function such that $\varphi'(t) > 0$,
$\varphi(t) - t$ is $2\pi$-periodic,
and any sequence $\{ \varphi(t + j\omega) \}$, $t \in \Rset$,
is a $p/q$-periodic configuration of $f$.
Then Proposition~\ref{prop:PeriodicConfiguration} implies
that the normal function $\psi(t)$ and
the incidence-reflection function $\theta(t)$
obtained from $\varphi(t)$ by relations~\eqref{eq:Relations}
satisfy the difference equation~\eqref{eq:DynamicEquation}.

\item[{\rm ({\bf p})}]
We define
\[
\lambda(t) =
 -\partial_1 S(\varphi(t),\varphi(t+\omega)) =
 \partial_2 S(\varphi(t-\omega),\varphi(t)).
\]
Since $\varphi(t) - t$ and $\lambda(t)$ are $2\pi$-periodic,
the map $c = (\varphi,\lambda) : \Rset \to \Rset^2$
can be projected to a map from $\Tset$ to $\Aset$.
Since $g(\varphi)$ is the support function of the caustic,
we get that $\lambda(t) = g(\varphi(t))$ and
$c(\Tset) = G \subset \Aset_\Gamma$.
Condition $\varphi'(t) > 0$ implies that $c'(t) \neq (0,0)$,
so $c:\Tset \to G$ is a parametrization.
Implicit equations~\eqref{eq:ImplicitEquations} imply
$f(c(t)) = c(t+\omega)$.

\item[{\rm ({\bf a})}]
The $p/q$-periodic sequences
$\{ \varphi(t + j\omega) \}_{j \in \Zset}$
form a $1$-parameter family of critical points of the
$p/q$-periodic action,
being $t \in \Rset$ a smooth parameter.
Therefore, the action is constant on this $1$-parameter family.
That is, all periodic billiard trajectories tangent
to the $p/q$-resonant caustic have the same length.
\qedhere
\end{enumerate}
\end{proof}

\begin{remark}
\label{rem:OneFormula}
If these three necessary conditions hold, then
$\mu \{ h' \circ \psi \cdot \sin \theta \} = 0$.
This relation follows by applying operator $\mu$
to~\eqref{eq:DynamicEquation}, since $\mu \circ \sigma = 2\mu$
and $\mu \circ \delta = 0$ on the space of
smooth $2\pi$-periodic functions.
\end{remark}

To provide a first insight into the usefulness of
Condition~({\bf f}) in Theorem~\ref{thm:ExistenceConditions},
let us give some information about functions
$\varphi(t)$, $\psi(t)$ and $\theta(t)$ in five examples.

\begin{example}
\label{exa:Circular}
The simplest example is the completely integrable circular billiard.
If $\Gamma$ is a circle of radius one centered at the origin,
then $h(\psi) \equiv 1$, so $\Aset_\Gamma = \Tset \times (-1,1)$
and the billiard map $f:\Aset_\Gamma \to \Aset_\Gamma$ is given by
$f(\varphi,\lambda) = (\varphi + \varpi(\lambda),\lambda)$
with $\varpi(\lambda) = 2 \arccos \lambda$.
In particular,
we can take $\varphi(t) = t - \omega/2$, $\psi(t) = t$,
$\theta(t) \equiv \omega/2$, $\omega = 2 \pi p/q$
and $L = 2 q \sin(\omega/2)$ for any $p/q \in (0,1/2)$.
Straightforward computations show that these functions
satisfy Condition~({\bf f}).
\end{example}

\begin{example}
Elliptic billiards are integrable too,
but their computations are harder.
If $\Gamma$ is the ellipse $\{ x^2/a^2 + y^2/b^2 = 1\}$,
then
$h(\psi) = \sqrt{a^2 \cos^2 \psi + b^2 \sin^2 \psi}$;
see~\cite[Lemma~1]{Bialy2022}.
The explicit expression of
$\varphi(t)$, $\psi(t)$ and $\theta(t)$ as functions of
a dynamical parameter $t$ requires the use of elliptic functions
whose modulus depends on the eccentricity of the ellipse
and the rotation number $p/q$ of each resonant caustic.
See~\cite{PintoRamirezRos2013,Damasceno_etal2017,Bialy2022}
for similar computations.
We omit the details,
since we only deal with deformations of circles in this work,
but we stress that $\varphi'(t),\psi'(t) \not \equiv 1$
in elliptic billiards, unlike circular billiards.
\end{example}

\begin{example}
\label{exa:ConstantWidth}
Constant width curves are a classic example.
Any such curve, other than a circle, has a nonsmooth (with cusps)
and nonconvex $1/2$-resonant caustic~\cite{Knill1998}.
The curve $\Gamma$ has constant width $w > 0$ when
$h(\psi) + h(\psi + \pi) \equiv w$, which implies that
all $(2\Zset \setminus \{0 \})$-harmonics of $h(\psi)$ vanish,
In that case, we take
$\varphi(t) = t - \omega/2$, $\psi(t) = t$,
$\theta(t) \equiv \omega/2$, $\omega = \pi$ and $L = 2w$
for $p/q = 1/2$.
These functions satisfy Condition~({\bf f})
because $\sin \theta(t) \equiv 1$, $\cos \theta(t) \equiv 0$ and
$h'(\psi) + h'(\psi+\pi) \equiv 0$.
\end{example}

\begin{example}
\emph{Gutkin billiard tables}~\cite{Gutkin2012,Bialy2018},
also called \emph{constant angle curves},
are another classic example.
We claim that circles are the only convex billiard tables
with a $p/q$-resonant caustic, with $p/q \in (0,1/2)$,
whose incidence-reflection function $\theta(t)$ is constant.
In that case, $\theta(t) \equiv \omega/2$ with $\omega = 2\pi p/q$,
so $\varphi(t) = t - \omega/2$ and $\psi(t) = t$
are the functions determined from $\theta(t)$
by relations~\eqref{eq:Relations}.
Therefore, the difference equation~\eqref{eq:DynamicEquation} becomes
\[
\tan(\omega/2) \big( h'(t + \omega) + h'(t) \big) =
h(t + \omega) - h(t).
\]
If $h(t) = \sum_{l \in \Zset} \hat{h}_l \rme^{\rmi l t}$
satisfies this equation,
then $\hat{h}_l = 0$ for any index $l \in \Zset$ such that
\[
\tan(l \pi p/q) \neq l \tan (\pi p /q).
\]
Cyr~\cite{Cyr2012} proved that given any integer $l \not \in \{-1,0,1\}$,
equation $\tan(l \pi p/q) = l \tan (\pi p /q)$
has no rational solution $p/q \in (0,1/2)$.
This proves the claim, because circles are the only convex curves
whose support function is a trigonometric polynomial
of degree one.
\end{example}

\begin{example}
A smooth centrally symmetric convex curve $\Gamma$
with support function $h(\psi)$
has a convex $1/4$-resonant caustic if and only if
\[
h^2(\psi) =
\hat{c}_0 + \sum_{l \in 2 + 4\Zset} \hat{c}_l \rme^{\rmi l \psi}
\qquad \mbox{ and } \qquad
h + h'' > 0.
\]
This claim is proved in~\cite[Proposition~3.1]{BialyTsodikovich2023a}.
Along the proof, the authors check three facts.
First, all such curves are centrally symmetric:
$h(\psi + \pi) = h(\psi)$.
Second, once fixed one of such curves, 
there is a constant $R > 0$ such that
$h^2(\psi) + h^2(\psi + \pi/2) = R^2$.
Third, then we can take $\psi(t) = t$ and determine
the incidence-reflection function $\theta(t)$
by means of
\[
h(t) = R \sin \theta(t), \quad
h(t+\pi/2) = R \cos \theta(t).
\]
Note that $\theta(t + \pi/2) = \pi/2 - \theta(t)$ and
$\psi(t + \pi/2) = \psi(t) + \pi/2$.
The last relation means that the tangent lines to $\Gamma$
at the impacts of any $1/4$-periodic trajectory form a rectangle.
This fact plays a key role in the Bialy-Mironov proof of a
strong version of the Birkhoff conjecture~\cite{BialyMironov2022}.
If we set $\varphi = \psi - \theta$ and $\omega = \pi/2$,
then functions $\varphi(t)$, $\psi(t)$ and $\theta(t)$
satisfy Condition~({\bf f}).
\end{example}

\section{Persistence of smooth convex resonant caustics}
\label{sec:Persistence}

Necessary conditions in Theorem~\ref{thm:ExistenceConditions}
are not sufficient for the existence of \emph{smooth convex}
resonant caustics, see Example~\ref{exa:ConstantWidth}.
However,
if the envelope of the 1-parameter family of lines from
$z(\psi(t))$ to $z(\psi(t+\omega))$,
where $t$ is the dynamical parameter,
is a smooth convex curve, then those conditions are sufficient too.
That is the case when we consider small enough smooth deformations
$\Gamma_\epsilon = \Gamma_0 + \Order(\epsilon)$ of a
smooth strictly convex curve $\Gamma_0$, not necessarily a circle,
with a smooth strictly convex resonant caustic.

We consider that setting.
To be precise, we assume the following hypothesis from here on.

\begin{enumerate}
\item[{\rm ({\bf H})}]
Let $p/q \in (0,1/2)$ be a rational rotational number such that
$\gcd(p,q) = 1$.
Let $\Gamma_0$ be a smooth strictly convex curve with
support function $h_0(\varphi)$.
We assume that there is a smooth convex $p/q$-resonant caustic
with support function $g_0(\varphi)$ inside $\Gamma_0$.
We also assume that the origin is in
the interior of the caustic, so $0 < g_0(\varphi) < h_0(\varphi)$.
Let $\Gamma_\epsilon = \Gamma_0 + \Order(\epsilon)$,
with $\epsilon \in [-\epsilon_0,\epsilon_0]$,
be a deformation of the unperturbed curve with
smooth support function $h(\varphi;\epsilon)$.
Let $(\varphi_1,\lambda_1) = f_\epsilon(\varphi,\lambda)$,
$S_\epsilon(\varphi,\varphi_1) = 2 h(\psi;\epsilon) \sin \theta$ and
$A_\epsilon\{\varphi\} =
 q \mu \{ S_\epsilon(\varphi,\tau\{\varphi\}) \} =
 2 q \mu \{ h \circ \psi \cdot \sin \theta \}$
be the perturbed billiard map in $\Gamma_\epsilon$,
the perturbed generating function and the perturbed action,
respectively.
\end{enumerate}

We need two parameters in that perturbed setting.
The \emph{dynamical parameter} $t$ parametrizes the invariant objects.
The \emph{perturbative parameter}
$\epsilon \in [-\epsilon_0,\epsilon_0]$ labels the ovals.
Then the shift, sum, difference and average operators are applied
to functions that depend on $t$ and $\epsilon$,
although they only act on $t$.
For instance,
$\tau\{ a(t;\epsilon) \} = a(t + \omega;\epsilon)$.
We will denote the derivatives of the support function
$h(\psi;\epsilon) = h_\epsilon(\psi)$ as
$h' = \frac{\rmd h}{\rmd \psi}$ and
$\dot{h} = \frac{\rmd h}{\rmd \epsilon}$.
Analogously, we will denote the derivatives of any function
$a(t;\epsilon) = a_\epsilon(t)$ as $a' = \frac{\rmd a}{\rmd t}$ and
$\dot{a} = \frac{\rmd a}{\rmd \epsilon}$.

Next, we state an immediate extension of
Theorem~\ref{thm:ExistenceConditions}.

\begin{corollary}
\label{cor:Persistence}
If $0 < \epsilon_0 \ll 1$,
the unperturbed smooth convex resonant $p/q$-caustic
persists under deformation $\Gamma_\epsilon$,
$\epsilon \in [-\epsilon_0,\epsilon_0]$,
if and only if the following three conditions hold.
\begin{enumerate}
\item[{\rm ({\bf F})}]
There are three smooth $2\pi$-periodic functions
$\varphi(t;\epsilon) - t$, $\psi(t;\epsilon) - t$
and $\theta(t;\epsilon)$ related by~\eqref{eq:Relations}
that satisfy the difference equation~\eqref{eq:DynamicEquation}.
Besides, $\varphi'(t;\epsilon) > 0$.

\item[{\rm ({\bf P})}]
There are smooth parametrizations
$c_\epsilon:\Tset \to \Aset_{\Gamma_\epsilon}$ such that
$G_\epsilon = c_\epsilon(\Tset)$ are graphs and
\[
f_\epsilon (c_\epsilon(t)) = c_\epsilon(t+\omega).
\]

\item[{\rm ({\bf A})}]
The $p/q$-periodic action is constant on the side function:
$A_\epsilon \{\varphi(t;\epsilon)\} = L(\epsilon)$,
where $L(\epsilon)$ is the length of all
$p/q$-periodic billiard trajectories in $\Gamma_\epsilon$
for $\epsilon \in [-\epsilon_0,\epsilon_0]$.
\end{enumerate}
\end{corollary}

\begin{remark}
\label{rem:TwoFormula}
Similarly to Remark~\ref{rem:OneFormula},
if these necessary conditions hold,
then
\[
2q \mu\{ h \circ \psi \cdot \sin \theta \} = L,\quad
\mu \{ h' \circ \psi \cdot \sin \theta \} = 0, \quad
2q \mu\{ \dot{h} \circ \psi \cdot \sin \theta \} = \dot{L}.
\]
Only the last formula is new.
Let us prove it.
If we derive the first relation with respect to $\epsilon$,
use the summation by parts formula and
take advantage of~\eqref{eq:DynamicEquation},
we get
\begin{align*}
\dot{L}
&=
2q\mu\big \{
\dot{h} \circ \psi \cdot \sin \theta +
h'\circ \psi \cdot \sin \theta \cdot \dot{\psi} +
h \circ \psi \cdot \cos \theta \cdot \dot{\theta} \big\} \\
&=
2q\mu\{ \dot{h} \circ \psi \cdot \sin \theta \} +
q\mu\big\{
h'\circ \psi \cdot \sin \theta \cdot \sigma\{ \dot{\varphi} \} +
h \circ \psi \cdot \cos \theta \cdot \delta\{ \dot{\varphi} \} \big\} \\
&=
2q\mu\{ \dot{h} \circ \psi \cdot \sin \theta \} +
q\mu\big\{
\big[
\sigma\{ h'\circ \psi \cdot \sin \theta \} -
\delta\{ h\circ \psi \cdot \cos\theta \}
\big] \cdot \tau\{\dot{\varphi}\}
\big\} \\
&=
2q\mu\{ \dot{h} \circ \psi \cdot \sin \theta \}.
\end{align*}
We can normalize $\Gamma_\epsilon$ by a scaling
in such a way that $\dot{L} = 0$, but we do not need it.
\end{remark}

Full persistence of resonant caustics is a extremely rare phenomenon,
so we introduce the more common concept of
\emph{$\Order(\epsilon^m)$-persistence}
for some \emph{order} $m \in \Nset \cup \{ 0 \}$.

\begin{definition}
\label{defi:PersistenceOrder}
Let $m \in \Nset \cup \{ 0 \}$.
The unperturbed resonant $p/q$-caustic
\emph{$\Order(\epsilon^m)$-persists}
under deformation $\Gamma_\epsilon$
if and only if the following three conditions hold.
\begin{enumerate}
\item[{\rm ({\bf F})$_m$}]
There are smooth $2\pi$-periodic functions
$\varphi(t;\epsilon) - t$, $\psi(t;\epsilon) - t$
and $\theta(t;\epsilon)$ related by~\eqref{eq:Relations}
such that $\varphi'(t;\epsilon) > 0$ and
\begin{equation}
\label{eq:DynamicEquationOrder}
\sigma \{ h'\circ \psi \cdot \sin \theta \} -
\delta \{ h\circ \psi \cdot \cos \theta \} =
\Order(\epsilon^{m+1})
\quad \mbox{ as } \epsilon \to 0.
\end{equation}

\item[{\rm ({\bf P})$_m$}]
There are smooth parametrizations
$c_\epsilon:\Tset \to \Aset_{\Gamma_\epsilon}$ such that
$G_\epsilon = c_\epsilon(\Tset)$ are graphs and
\[
f_\epsilon \circ c_\epsilon - \tau\{ c_\epsilon \} =
 \big( 0,\Order(\epsilon^{m+1}) \big) \quad
\mbox{ as } \epsilon \to 0.
\]

\item[{\rm ({\bf A})$_m$}]
There is a `length' $L(\epsilon) > 0$ such that
the side function $\varphi = \varphi(t;\epsilon)$ satisfies
\[
A_\epsilon \{\varphi\} - L(\epsilon) =
\Order(\epsilon^{m+1})
\quad \mbox{ as } \epsilon \to 0.
\]
\end{enumerate}
If the resonant caustic $\Order(\epsilon^m)$-persists,
but not $\Order(\epsilon^{m+1})$-persists,
we say that deformation $\Gamma_\epsilon$
\emph{$\Order(\epsilon^{m+1})$-breaks} the caustic.
\end{definition}

Since $\Order(\epsilon^m)$-persistence is the main concept of this work,
some comments are in order.
If $\varphi_0(t)$, $\psi_0(t)$ and $\theta_0(t)$ are the
unperturbed side, normal and incidence-reflection functions, then
\[
\sigma \{ h'\circ \psi_0 \cdot \sin \theta_0 \} -
\delta \{ h\circ \psi_0 \cdot \cos \theta_0 \} =
\Order(\epsilon)
\quad \mbox{ as } \epsilon \to 0,
\]
since $h = h_0 + \Order(\epsilon)$.
Therefore, 
the unperturbed caustic always $\Order(\epsilon^0)$-persists.
We do not need to check \emph{all} three Conditions~({\bf F})$_m$,
({\bf P})$_m$ and~({\bf A})$_m$,
because there are logical dependencies among them.
We prove in Proposition~\ref{prop:ErrorMelnikov} that
Conditions~({\bf F})$_m$ and~({\bf P})$_m$ are equivalent and
both imply Condition~({\bf A})$_m$.
We will only check Condition~({\bf F})$_m$ in our computations.
We have included the other conditions as part of our definition
to present a broader view of the problem.

The three conditions look similar,
but they have different characteristics.
On the one hand,
Conditions~({\bf F})$_m$ and~({\bf P})$_m$
are stated in terms of a \emph{single} iteration of
the perturbed map $f_\epsilon$.
On the other hand,
Condition~({\bf A})$_m$ requires to consider all the shifts
\[
\bar{\varphi}_j =
\bar{\varphi}_j(t;\epsilon) =
\tau^j\{ \varphi(t;\epsilon) \} =
\varphi(t + j\omega;\epsilon), \quad j=0,\ldots,q.
\]
Hence, Conditions~({\bf F})$_m$ and~({\bf P})$_m$ are
easier to deal with from a computational point of view.

Condition~({\bf F})$_m$ means that there is a
reparametrization $\varphi_\epsilon(t) = \varphi(t;\epsilon)$
of the original angle $\varphi$
in terms of a new dynamical parameter $t$ such that
\[
\partial_2 S_\epsilon(\varphi_\epsilon(t-\omega),\varphi_\epsilon(t)) +
\partial_1 S_\epsilon(\varphi_\epsilon(t),\varphi_\epsilon(t+\omega)) =
\Order(\epsilon^{m+1})
\quad \mbox{ as } \epsilon \to 0.
\]
Thus, Condition~({\bf F})$_m$ is related to
the Lagrangian formulation as a second-order difference equation
of the invariance condition for (nonresonant and resonant)
rotational curves of exact twist maps.
It was inspired by ideas contained
in~\cite{LeviMoser2001,KaloshinZhang2018}.
Compare with~\cite[Eq. (3)]{LeviMoser2001} in the nonresonant setting and
with \cite[Eq. (1)]{KaloshinZhang2018} in the resonant setting.
Finally, Condition~({\bf A})$_m$ follows the variational
approach in~\cite{Damasceno_etal2017}.

Let us prove the logical dependencies among these
$\Order(\epsilon^m)$-persistence conditions
and how the dominant terms in their
$\Order(\epsilon^{m+1})$-errors are related.

\begin{proposition}
\label{prop:ErrorMelnikov}
Conditions~{\rm ({\bf F})$_m$} and~{\rm ({\bf P})$_m$} are equivalent
and both imply Condition~{\rm ({\bf A})$_m$}.
Let $\varphi_0(t) = \varphi(t;0)$ be
the unperturbed side function.
If the $p/q$-resonant caustic \emph{$\Order(\epsilon^m)$-persists}
and we follow the notations introduced in
Definition~\ref{defi:PersistenceOrder},
then there is a smooth $2\pi$-periodic function
$\mathcal{E}_{m+1}(t)$ and a smooth $2\pi/q$-periodic function
$\mathcal{L}_{m+1}(t)$ such that
\begin{subequations}
\label{eq:seq}
\label{equations}
\begin{align}
\label{seq:F}
\sigma \{ h'\circ \psi \cdot \sin \theta \} -
\delta \{ h\circ \psi \cdot \cos \theta \} &=
\epsilon^{m+1} \tau\{\mathcal{E}_{m+1}\} + \Order(\epsilon^{m+2}), \\
\label{seq:P}
f_\epsilon \circ c_\epsilon - \tau\{ c_\epsilon \} &=
\big(
0, \epsilon^{m+1} \tau\{\mathcal{E}_{m+1}\} + \Order(\epsilon^{m+2})
\big), \\
\label{seq:A}
A_\epsilon \{ \varphi \} - L(\epsilon) &=
\epsilon^{m+1} q \mathcal{L}_{m+1} + \Order(\epsilon^{m+2}),
\end{align}
\end{subequations}
as $\epsilon \to 0$.
Besides,
$\mu\{ \mathcal{E}_{m+1} \varphi'_0 \} = \mathcal{L}'_{m+1}$ and
$\int_\Tset \mathcal{E}_{m+1} \varphi'_0 = 0$.
\end{proposition}

\begin{proof}
Firstly, we check that ({\bf F})$_m \Rightarrow$
({\bf P})$_m$ \& ({\bf A})$_m$,
$\mathcal{L}'_{m+1} = \mu\{ \mathcal{E}_{m+1} \varphi'_0 \}$ and
$\int_\Tset \mathcal{E}_{m+1} \varphi'_0 = 0$.

Let $\varphi(t;\epsilon)$, $\psi(t;\epsilon)$ and $\theta(t;\epsilon)$
be the functions described in ({\bf F})$_m$ and $\mathcal{E}_{m+1}(t)$
be the function determined by~\eqref{seq:F}.
Set $\varphi_1(t;\epsilon) = \varphi(t+\omega;\epsilon)$,
so $\varphi_1 = \tau\{ \varphi \}$.
We consider the smooth $2\pi$-periodic functions
$\lambda(t;\epsilon)$ and $\lambda_1(t;\epsilon)$ given by
\[
\lambda =
h \circ \psi \cdot \cos \theta - h' \circ \psi \cdot \sin \theta, \qquad
\lambda_1 =
h \circ \psi \cdot \cos \theta + h' \circ \psi \cdot \sin \theta.
\]
Set
$c_\epsilon(t) =
 \big( \varphi(t;\epsilon),\lambda(t;\epsilon) \big)$ and
$d_\epsilon(t) =
 \big( \varphi_1(t;\epsilon),\lambda_1(t;\epsilon) \big)$.
Implicit equations~\eqref{eq:ImplicitEquations} imply that
$f \circ c_\epsilon = d_\epsilon$.
Besides,
$\lambda_1 - \tau\{ \lambda \} =
\sigma\{ h' \circ \psi \cdot \sin \theta \} -
\delta\{ h \circ \psi \cdot \cos \theta \} =
\epsilon^{m+1} \tau\{ \mathcal{E}_{m+1} \} + \Order(\epsilon^{m+2})$,
which is equivalent to estimate~\eqref{seq:P}.
This proves ({\bf P})$_m$.

Set $\bar{\varphi}_j = \tau^j \{ \varphi \}$,
$\bar{\psi}_j = \tau^j \{ \psi \}$ and
$\bar{\theta}_j = \tau^j \{ \theta \}$ for all $j \in \Zset$.
Note that
$\bar{\varphi}_{j+q}(t;\epsilon) =
 \bar{\varphi}_j(t;\epsilon) + 2\pi p$,
$\bar{\psi}_{j+q}(t;\epsilon) =
 \bar{\psi}_j(t;\epsilon) + 2\pi p$ and
$\bar{\theta}_{j+q}(t;\epsilon) = \bar{\theta}_j(t;\epsilon)$.
If we derive the expression that defines the action and
we recall that $\varphi = \varphi_0 + \Order(\epsilon)$,
then we obtain the estimate
\begin{align*}
\frac{\rmd}{\rmd t}\big[ A_\epsilon \{ \varphi \} \big]
&=
\textstyle \sum_{j=0}^{q-1}
\big[
\partial_1 S_\epsilon (\bar{\varphi}_j,\bar{\varphi}_{j+1})
\bar{\varphi}'_j +
\partial_2 S_\epsilon (\bar{\varphi}_j,\bar{\varphi}_{j+1})
\bar{\varphi}'_{j+1}
\big] \\
&=
\textstyle \sum_{j=0}^{q-1}
\big[
\partial_2 S_\epsilon (\bar{\varphi}_{j-1},\bar{\varphi}_j) +
\partial_1 S_\epsilon (\bar{\varphi}_j,\bar{\varphi}_{j+1})
\big]\bar{\varphi}'_j \\
&=
\textstyle \sum_{j=0}^{q-1}
\big[
\sigma \{ h'\circ \bar{\psi}_{j-1} \cdot \sin \bar{\theta}_{j-1} \} -
\delta \{ h\circ \bar{\psi}_{j-1} \cdot \cos \bar{\theta}_{j-1} \}
\big] \bar{\varphi}'_j \\
&=
\epsilon^{m+1}
\textstyle \sum_{j=0}^{q-1} \tau^j\{ \mathcal{E}_{m+1} \varphi' \} +
\Order(\epsilon^{m+2}) \\
&=
\epsilon^{m+1} q \mu\{ \mathcal{E}_{m+1}  \varphi'\} +
\Order(\epsilon^{m+2}) \\
&=
\epsilon^{m+1} q \mu\{ \mathcal{E}_{m+1}  \varphi'_0 \} +
\Order(\epsilon^{m+2}),
\end{align*}
which, by integration, is equivalent to estimate~\eqref{seq:A}
for any smooth $2\pi/q$-periodic function $\mathcal{L}_{m+1}$
such that $\mathcal{L}'_{m+1} = \mu\{ \mathcal{E}_{m+1} \varphi'_0 \}$.
This proves~({\bf A})$_m$ and the relation between
$\mathcal{L}_{m+1}$ and $\mathcal{E}_{m+1}$.
The operator $\mu$ is the projection onto the $q\Zset$-harmonics,
so the zero-th harmonics of $\mathcal{E}_{m+1} \varphi'_0$
and $\mathcal{L}'_{m+1}$ coincide, so
$\int_\Tset \mathcal{E}_{m+1} \varphi'_0 =
 \int_\Tset \mathcal{L}'_{m+1} = 0$.

Secondly, we check that ({\bf P})$_m \Rightarrow$ ({\bf F})$_m$.
Let
$c_\epsilon(t) = \big( \varphi(t;\epsilon),\lambda(t;\epsilon) \big)$
be the parametrization described in ({\bf P})$_m$ and
$\mathcal{E}_{m+1}(t)$ be the error function given in~\eqref{seq:P}.
Property $\varphi' > 0$ holds because
$c_\epsilon:
 \Tset \to G_\epsilon \subset
           \Aset_{\Gamma_\epsilon}$ is a
parametrization and $G_\epsilon$ is a graph.
Functions $\varphi$, $\psi = \sigma\{\varphi\}/2$ and
$\theta = \delta\{ \varphi \}/2$ satisfy relations~\eqref{eq:Relations}.
Set
$\big( \varphi_1(t;\epsilon),\lambda_1(t;\epsilon) \big) =
 f_\epsilon(c_\epsilon(t))$.
We deduce from implicit equations~\eqref{eq:ImplicitEquations} that
\[
\sigma \{ h'\circ \psi \cdot \sin \theta \} -
\delta \{ h\circ \psi \cdot \cos \theta \} =
\lambda_1 - \tau\{ \lambda \} =
\epsilon^{m+1} \tau\{ \mathcal{E}_{m+1} \} + \Order(\epsilon^{m+2}),
\]
which is exactly estimate~\eqref{seq:F}.
This proves~({\bf F})$_m$.
\end{proof}

\begin{remark}
\label{rem:Flux}
We can prove that $\int_\Tset \mathcal{E}_{m+1} \varphi'_0 = 0$
in another way.
The flux of the exact twist map $f_\epsilon$ across
the graph $G_\epsilon = c_\epsilon(\Tset)$,
which is a rotational curve, is zero~\cite[\S V]{Meiss1992}.
This flux is
\begin{align*}
\int_{f_\epsilon(G_\epsilon)} \lambda \rmd \varphi -
\int_{G_\epsilon} \lambda \rmd \varphi
&=
\int_\Tset (\lambda_1 \varphi'_1 - \lambda \varphi') \\
&=
\int_\Tset \delta \{ \lambda \varphi'\} +
\epsilon^{m+1} \int_\Tset \tau\{ \mathcal{E}_{m+1} \varphi' \} +
\Order(\epsilon^{m+2}) \\
&=
\epsilon^{m+1} \int_\Tset \mathcal{E}_{m+1} \varphi'_0 +
\Order(\epsilon^{m+2}) \quad
\mbox{ as } \epsilon \to 0.
\end{align*}
We have used that if
$c_\epsilon(t) =
 \big( \varphi(t;\epsilon),\lambda(t;\epsilon) \big)$ and
$(\varphi_1,\lambda_1) = f(\varphi,\lambda)$,
then $\varphi_1 = \tau\{ \varphi \}$ and
$\lambda_1 = \tau\{ \lambda \} +
\epsilon^{m+1} \tau\{ \mathcal{E}_{m+1} \} + \Order(\epsilon^{m+2})$.
We have also used that $\varphi = \varphi_0 + \Order(\epsilon)$.
\end{remark}

\begin{remark}
\label{rem:Graphs}
Following~\cite{RamirezRos2006,PintoRamirezRos2013},
we could also have considered a fourth condition defined in terms
of $\varphi \in \Tset$ instead of $t \in \Rset$.
However, such approach forces us to deal with
the power map $f^q_\epsilon$, which is technically impractical.
We have not pursued it.
That discarded condition is:
\begin{enumerate}
\item[{\rm ({\bf G})$_m$}]
There are smooth functions
$g^\bullet_\epsilon,g^\star_\epsilon:\Tset \to \Rset$ such that
$f^q_\epsilon(\varphi,g^\bullet_\epsilon(\varphi)) =
(\varphi,g^\star_\epsilon(\varphi))$,
$g^\star_\epsilon, g^\bullet_\epsilon = g_0 + \Order(\epsilon)$ and
$g^\star_\epsilon - g^\bullet_\epsilon = \Order(\epsilon^{m+1})$
as $\epsilon \to 0$.
\end{enumerate}
It means that $f^q_\epsilon$ projects
$G^\bullet_\epsilon = \Graph(g^\bullet_\epsilon)$ onto
$G^\star_\epsilon = \Graph(g^\star_\epsilon)$ in
the vertical direction.
None of these two graphs have to coincide with
$G_\epsilon = c_\epsilon(\Tset)$,
but all of them are $\Order(\epsilon^{m+1})$-close.
Conditions~({\bf F})$_m$ or~({\bf P})$_m$
imply Condition~({\bf G})$_m$.
The proof of this implication is based on two ideas.
First, we check that the restriction of the power map
$f^q_\epsilon$ to the graph $G_\epsilon$
is $\Order(\epsilon^{m+1})$-close to the identity,
although the error in the horizontal direction may be nonzero.
Next, we use the twist condition to cancel the horizontal error,
so we get the graphs
$G^\bullet_\epsilon, G^\star_\epsilon =
 G_\epsilon + \Order(\epsilon^{m+1})$.
We omit the technical details since we do not need this result.
On the contrary,
we do not know whether Condition~({\bf G})$_m$
implies Conditions~({\bf F})$_m$ and~({\bf P})$_m$,
because it is not clear how to extract information about
a single iterate of the perturbed map from information about
the power map $f^q_\epsilon$.
\end{remark}

\begin{definition}
\label{defi:Melnikov}
The $2\pi/q$-periodic function $\mathcal{L}_{m+1}(t)$ is the
\emph{$p/q$-resonant potential of order $m+1$}.
The $2\pi$-periodic function $\mathcal{E}_{m+1}(t)$ is the
\emph{$p/q$-resonant error of order $m+1$}.
\end{definition}

We will check in the next section that
if the $p/q$-resonant caustic $\Order(\epsilon^m)$-persists
under a deformation of the unit circle,
then the Melnikov potential $\mathcal{L}_{m+1}(t)$
is uniquely determined from previously
computed objects.
See the second item in Proposition \ref{prop:Iterative}.

\section{High-order persistence in deformed circles}
\label{sec:Circles}

Let us apply the previous high-order persistence
theory to smooth deformations of circles.
The main goal is to check that
the smooth $2\pi$-periodic coefficients of the Taylor
expansions in powers of the perturbative parameter $\epsilon$
of the perturbed side, normal and incidence-reflection functions
can be computed recursively order by order as long as some
compatibility conditions hold.
Such compatibility conditions have to do with the inversion
of the difference operator $\delta$,
so they boil down to the fact that certain smooth $2\pi$-periodic
functions have no $q\Zset$-resonant harmonics.
We look for a practical way to find the exact order
at which a given resonant caustic is destroyed, so we
write down explicit formulas for all recursive computations.

To begin with,
we assume the following hypothesis from here on.

\begin{enumerate}
\item[{\rm ({\bf H'})}]
Let $p/q \in (0,1/2)$ be a rational rotational number such that
$\gcd(p,q) = 1$.
Let $\Gamma_0$ be the unit circle centered at the origin,
see Example~\ref{exa:Circular}.
The unperturbed side, normal and incidence-reflection functions are
\[
\varphi_0(t) = t - \omega/2, \quad
\psi_0(t) = t, \quad
\theta_0(t) = \omega/2, \quad
\omega = 2\pi p/q.
\]
Let $\Gamma_\epsilon$,
with $\epsilon \in [-\epsilon_0,\epsilon_0]$,
be a deformation of $\Gamma_0$ with
smooth support function~\eqref{eq:SupportFunction}.
\end{enumerate}

We look for some perturbed side, normal and
incidence-reflection functions
\begin{equation}
\label{eq:FunctionsExpasions}
\varphi(t;\epsilon) \asymp 
\sum_{k \ge 0} \epsilon^k \varphi_k(t), \quad
\psi(t;\epsilon) \asymp
\sum_{k \ge 0} \epsilon^k \psi_k(t), \quad
\theta(t;\epsilon) \asymp
\sum_{k \ge 0} \epsilon^k \theta_k(t)
\end{equation}
that satisfy Condition~({\bf F})$_m$
for an order $m \in \Nset \cup \{0\}$ as high as possible.

\begin{nota}
If $a(t;\epsilon) \asymp \sum_{k \ge 0} a_k(t) \epsilon^k$
as $\epsilon \to 0$,
then $a_{\le m}(t;\epsilon) = \sum_{k=0}^m a_k(t) \epsilon^k$.
Symbol $a_{\le m}(t;\epsilon)$ may be used even when there are no
coefficients $a_k(t)$ with $k > m$.
Symbol $a_{<m}(t;\epsilon)$ has a similar meaning.
An expression like
$\mathcal{C}_m = \mathcal{C}_m\llbracket a_{\le m},b_{<m} \rrbracket$
means that $\mathcal{C}_m$ is a smooth $2\pi$-periodic
function that can be written as a \emph{differential} expression
in some smooth $2\pi$-periodic functions
$a_k$ for $1 \le k \le m$ and $b_k$ for $1 \le k < m$.
The term differential means that the derivatives of functions
$a_k$ and $b_k$ can appear in those expressions.
\end{nota}

The study of the approximate difference
equation~\eqref{eq:DynamicEquationOrder} requires
to recursively compute the asymptotic expansions of
$\mathcal{R} := h \circ \psi \cdot \cos \theta$ and
$\mathcal{Q} := h' \circ \psi \cdot \sin \theta$ as $\epsilon \to 0$.

\begin{lemma}
\label{lem:QR}
Set $c = \cos(\omega/2)$ and $s = \sin(\omega/2)$.
The coefficients of the asymptotic expansions
\[
\mathcal{R} = h \circ \psi \cdot \cos \theta \asymp
c + \sum_{k \ge 1} \mathcal{R}_k \epsilon^k, \quad
\mathcal{Q} = h' \circ \psi \cdot \sin \theta \asymp
\sum_{k \ge 1} \mathcal{Q}_k \epsilon^k
\]
have the form
\begin{align*}
\mathcal{Q}_k &=
\mathcal{Q}_k \llbracket h'_{\le k},\psi_{<k},\theta_{<k} \rrbracket =
s h'_k +
\tilde{\mathcal{Q}}_k \llbracket h'_{<k},\psi_{<k},\theta_{<k} \rrbracket, \\
\mathcal{R}_k &=
 \mathcal{R}_k \llbracket h_{\le k},\psi_{<k},\theta_{\le k} \rrbracket =
\tilde{\mathcal{R}}_k \llbracket h_{\le k},\psi_{<k},\theta_{<k} \rrbracket -
s\theta_k.
\end{align*}
Besides,
$\tilde{\mathcal{Q}}_1 = 0$,
$\tilde{\mathcal{R}}_1 = ch_1$,
$\tilde{\mathcal{Q}}_2 = s h''_1 \psi_1 + c h'_1 \theta_1$ and
$\tilde{\mathcal{R}}_2 = c h_2 + c h'_1 \psi_1 - s h_1 \theta_1 - c \theta^2_1/2$.
\end{lemma}

The proof of Lemma~\ref{lem:QR} is postponed to Appendix~\ref{app:Expansions}.
The first coefficients are obtained from the
explicit recurrences for $\tilde{\mathcal{Q}}_k$ and $\tilde{\mathcal{R}}_k$
given in Lemma~\ref{lem:QRbis}.

Once we know that $\mathcal{Q}_{< k}$ and $\mathcal{R}_{< k}$
only depend on $\psi_{<k}$ and $\theta_{< k}$,
we deduce that if
$\varphi_{<k}(t;\epsilon)$, $\psi_{<k}(t;\epsilon)$ and
$\theta_{<k}(t;\epsilon)$ satisfy Condition ({\bf F})$_{k-1}$,
then
$\varphi_{\le k}(t;\epsilon)$, $\psi_{\le k}(t;\epsilon)$ and
$\theta_{\le k}(t;\epsilon)$ satisfy Condition ({\bf F})$_k$
provided coefficients
$\varphi_k(t)$, $\psi_k(t)$ and $\theta_k(t)$ are chosen
in such a way that
$\sigma\{ \mathcal{Q}_k \} = \delta\{ \mathcal{R}_k \}$.
We prove below that these $k$-th coefficients can be found
if and only if $\mathcal{Q}_k$ has no $q\Zset$-resonant harmonics,
in which case all three $k$-th coefficients are uniquely determined
provided that they have no $q\Zset$-resonant harmonics either.

\begin{proposition}
\label{prop:Iterative}
Let $k \in \Nset$ be a fixed order.
If the $p/q$-resonant caustic $\Order(\epsilon^{k-1})$-persists,
then the following properties hold.
\begin{enumerate}[a)]
\item
$\bar{\mathcal{Q}}_k := \frac{1}{2\pi} \int_\Tset \mathcal{Q}_k = 0$.

\item
The $p/q$-resonant Melnikov potential of order $k$
is completely determined from previously computed objects:
$\mathcal{L}_k =
 \mathcal{L}_k
\llbracket h_{\le k},\psi_{<k},\theta_{<k} \rrbracket$.

\item
The $p/q$-resonant caustic $\Order(\epsilon^k)$-persists
if and only if
\begin{equation}
\label{eq:PersistenceCondition}
\mu\{ \mathcal{Q}_k \} = 0,
\end{equation}
in which case $\theta_k(t)$ is
the unique smooth $2\pi$-periodic solution of
\begin{equation}
\label{eq:CompatibleEquations}
s \delta\{ \theta_k \} =
\delta\{ \tilde{\mathcal{R}}_k \} - \sigma\{ \mathcal{Q}_k \},\qquad
\mu\{ \theta_k \} = 0,
\end{equation}
and then $\varphi_k(t)$ and $\psi_k(t)$
are uniquely determined from $\theta_k(t)$
by
\begin{equation}
\label{eq:RelationsCoefficients}
\delta\{ \varphi_k \} = 2\theta_k, \qquad
\mu\{ \varphi_k\} = 0, \qquad
\psi_k = \varphi_k + \theta_k.
\end{equation}
\item
If condition~\eqref{eq:PersistenceCondition} fails,
then it is satisfied for any support function
$h^\star = h + \epsilon^k \eta_k$ such that
\begin{equation}
\label{eq:CorrectionEquation}
\mu\big\{ s \eta'_k + \mathcal{Q}_k \big\} = 0,
\end{equation}
so the $p/q$-resonant caustic $\Order(\epsilon^k)$-persists
under a corrected deformation
$\Gamma^\star_\epsilon = \Gamma_\epsilon + \Order(\epsilon^k)$.
We can choose the correction $\eta_k$ in such a way that it
only contains $(q\Zset \setminus \{0\})$-harmonics.
\end{enumerate}
\end{proposition}

\begin{proof}
We use several times  along this proof that $\mu \circ \tau = \mu$,
$\mu \circ \sigma = 2\mu$ and $\mu \circ \delta = 0$
on the space of smooth $2\pi$-periodic functions.
We also use that if $b(t)$ is a smooth $2\pi$-periodic function,
$\delta\{ a \} = b$ has a smooth $2\pi$-periodic
solution $a(t)$ if and only if $\mu\{ b \} = 0$,
in which case the solution is unique under the additional condition
$\mu\{ a \} = 0$. See Lemma~\ref{lem:Inversion}.

\begin{enumerate}[a)]
\item
The asymptotic expansions in Lemma~\ref{lem:QR} imply that
\[
\sigma \{ h'\circ \psi \cdot \sin \theta \} -
\delta \{ h\circ \psi \cdot \cos \theta \} =
\sum_{j = 1}^k
\big( \sigma\{ \mathcal{Q}_j \} - \delta\{ \mathcal{R}_j \} \big)
\epsilon^j + \Order(\epsilon^{k+1}).
\]
The $\Order(\epsilon^{k-1})$-persistence hypothesis means
that $\sigma\{ \mathcal{Q}_j \} = \delta\{ \mathcal{R}_j \}$
for $j=1,\ldots,k-1$.
If $\mathcal{E}_k$ and $\mathcal{L}_k$ are the $p/q$-resonant
error and potential of order $k$, respectively, then
$\tau\{ \mathcal{E}_k \} = 
\sigma\{ \mathcal{Q}_k \} - \delta\{ \mathcal{R}_k \}$ and
$\mathcal{L}'_k =
 \mu\{ \mathcal{E}_k \} =
 2\mu \{ \mathcal{Q}_k \}$,
since $\varphi'_0(t) = 1$.
Finally,
$\bar{\mathcal{Q}}_k =
 \frac{1}{2\pi} \int_\Tset \mathcal{Q}_k =
 \frac{1}{4\pi} \int_\Tset \mathcal{E}_k = 0$.
(See the last claims in Proposition~\ref{prop:ErrorMelnikov}.)

\item
It is a direct consequence of properties
$\mathcal{L}'_k = 2\mu\{ \mathcal{Q}_k \}$ and
$\mathcal{Q}_k =
\mathcal{Q}_k
\llbracket h'_{\le k},\psi_{<k},\theta_{<k} \rrbracket$.
 
\item
Condition~\eqref{eq:PersistenceCondition} is necessary:
$\Order(\epsilon^k)$-persistence means
$\sigma\{ \mathcal{Q}_k \} = \delta\{ \mathcal{R}_k \}$,
so $\mu\{ \mathcal{Q}_k \} = 0$.
Condition~\eqref{eq:PersistenceCondition} is sufficient:
If $\mu\{ \mathcal{Q}_k \} = 0$, then
\[
\mu
\big \{
\delta\{ \tilde{\mathcal{R}}_k \} - \sigma\{ \mathcal{Q}_k \}
\big \} =
-2 \mu \{ \mathcal{Q}_k \} = 0,
\]
so~\eqref{eq:CompatibleEquations} has a unique
smooth $2\pi$-periodic solution $\theta_k(t)$.
Then we find $\varphi_k$ and $\psi_k$
by means of relations~\eqref{eq:RelationsCoefficients}.
Finally,
the $\Order(\epsilon^k)$-terms of
$\sigma \{ h'\circ \psi \cdot \sin \theta \}$ and
$\delta \{ h\circ \psi \cdot \cos \theta \}$ are
$\sigma\{ \mathcal{Q}_k \}$ and
$\delta\{ \tilde{\mathcal{R}}_k - s \theta_k \}$,
respectively.
See Lemma~\ref{lem:QR}.
Both terms coincide when the first
relation in~\eqref{eq:CompatibleEquations} holds.
Property $\varphi'(t;\epsilon) = 1 + \Order(\epsilon) > 0$ holds
for $\epsilon \in [-\epsilon_0,\epsilon_0]$
if $0 < \epsilon_0 \ll 1$.
Relations~\eqref{eq:Relations} follow from condition
$\mu\{ \theta_k \} = 0$ and relations~\eqref{eq:RelationsCoefficients},
since $\psi_k = \varphi_k + \theta_k$ implies that
$\mu \{ \psi_k \} = 0$ and $\psi_k = \sigma\{ \varphi_k \}/2$.
This proves that Condition~({\bf F})$_k$ holds.

\item
We deduce from Lemma~\ref{lem:QR} that
if $\mathcal{Q}^\star_k$ is the $\Order(\epsilon^k)$-coefficient
of $(h^\star)' \circ \psi_{<k} \cdot \sin \theta_{<k}$,
then $\mathcal{Q}^\star_k = s \eta'_k + \mathcal{Q}_k$.
The existence of the correction $\eta_k$ follows from
property $\bar{\mathcal{Q}}_k = 0$.
\qedhere
\end{enumerate}
\end{proof}

Finally, we prove Theorem~\ref{thm:Persistence} by recursively
applying Proposition~\ref{prop:Iterative} for $k=1,\ldots,m$.

\begin{proof}[Proof of Theorem~\ref{thm:Persistence}]
We denote the Fourier coefficients of smooth $2\pi$-periodic
functions with a hat, so we write
$h_1(t) =
 \sum_{l \in \Zset} \hat{h}_{1,l} \rme^{\rmi l t}$ and
$\theta_1(t) =
 \sum_{l \in \Zset} \hat{\theta}_{1,l} \rme^{\rmi l t}$.

\begin{enumerate}[a)]
\item
We know that $\mathcal{Q}_1 = s h'_1$ from Lemma~\ref{lem:QR}.
Hence, the necessary and sufficient condition~\eqref{eq:PersistenceCondition}
for $\Order(\epsilon)$-persistence becomes $\mu \{ h'_1 \} = 0$
or, equivalently, $\mu_{q\Zset^\star} \{ h_1 \} = 0$.

\item
We know that $\mathcal{Q}_2 = s h'_2 + \tilde{\mathcal{Q}}_2$
and $\tilde{\mathcal{Q}}_2 = s h''_1 \psi_1 + c h'_1 \theta_1$
from Lemma~\ref{lem:QR}.
Let us check that
$\mu\{ \tilde{\mathcal{Q}}_2 \} = s \mu\{ \theta_1 \theta'_1 \}$.
In order to do that, we need three properties.

Firstly,
once we know that the resonant caustic $\Order(\epsilon)$-persists,
we determine the first-order coefficient $\theta_1(t)$  by
solving~\eqref{eq:CompatibleEquations} with $k=1$.
That is,
\begin{equation}
\label{eq:CompatibleEquation1}
s \delta\{ \theta_1 \} =
\delta\{ \tilde{\mathcal{R}}_1 \} - \sigma\{ \mathcal{Q}_1 \} =
\delta\{ ch_1 \} - \sigma\{ sh'_1 \}, \qquad
\mu\{ \theta_1 \}  = 0.
\end{equation}
Secondly,
we determine $\varphi_1(t)$ and $\psi_1(t)$ by
solving~\eqref{eq:RelationsCoefficients} with $k=1$.
In particular,
\[
\theta_1 = \delta\{ \varphi_1 \}/2, \qquad
\psi_1 = \sigma\{ \varphi_1 \}/2.
\]
Thirdly,
we recall that if $a(t)$ and $b(t)$ are $2\pi$-periodic functions,
then
\[
\mu \big\{ a \delta\{b\} \big\} =
-\mu \big\{ \delta\{ a \} \tau \{ b \} \big\}, \qquad
\mu \big\{ a \sigma\{ b \} \big\} =
\mu \big\{ \sigma\{ a \} \tau\{ b \} \big\}.
\]
Next, we use these three properties to get the formula we are looking for:
\begin{align*}
\mu\{ \tilde{\mathcal{Q}}_2 \}
&=
\mu\{ s h''_1 \psi_1 + c h'_1 \theta_1 \} =
\mu\big\{
s h''_1 \sigma\{ \varphi_1\} + c h'_1 \delta\{ \varphi_1 \}
\big\}/2 \\
&=
\mu\big\{
\big[\sigma\{ s h''_1\} - \delta\{ c h'_1\} \big] \tau\{ \varphi_1 \}
\big\}/2 =
\mu\big\{
\big[\sigma\{ s h'_1\} - \delta\{ c h_1\} \big]' \tau\{ \varphi_1 \}
\big\}/2 \\
&=
-s \mu\big\{ \delta\{ \theta'_1\} \tau\{\varphi_1 \} \big\}/2 =
s \mu\big\{ \delta\{ \varphi_1\} \theta'_1 \big\}/2 =
s \mu\{ \theta_1 \theta'_1 \}.
\end{align*}
Thus, condition~\eqref{eq:PersistenceCondition}
for $\Order(\epsilon^2)$-persistence becomes
$\mu\big\{ h'_2 + \theta_1 \theta'_1 \big\} = 0$ or,
equivalently,
\[
\mu_{q \Zset^\star} \{ h_2 + \theta_1^2/2 \} = 0. 
\]
Set $e^\pm_l = \rme^{\rmi l \omega} \pm 1$.
The Fourier coefficients of the unique smooth $2\pi$-periodic
solution  of~\eqref{eq:CompatibleEquation1} are easily determined:
$\hat{\theta}_{1,l} = 0$ for all $l \in q\Zset$ and
\[
s e^-_{l} \hat{\theta}_{1,l} =
\big( c e^-_l - \rmi l s e^+_l \big) \hat{h}_{1,l}, \qquad
\forall l \not \in q\Zset.
\]
The RHS of this identity vanishes, and $\hat{\theta}_{1,l}$ too,
when $l = \pm 1$.
If $l \not \in q\Zset \cup \{-1,1\}$, then
\[
\hat{\theta}_{1,l} = 
\left( \frac{c}{s} -  \rmi l \frac{e^+_l}{e^-_l} \right)
\hat{h}_{1,l} =
\begin{cases}
\displaystyle \frac{\tan(l\omega/2) - l \tan(\omega/2)}
     {\tan(\omega/2) \tan(l\omega/2)}
\hat{h}_{1,l}, & \mbox{ if $2l \not \in q\Zset$}, \\ 
\displaystyle \frac{\hat{h}_{1,l}}{\tan(\omega/2)}, &
\mbox{ if $2l \in q\Zset$}.
\end{cases}
\]
We have used that $e^+_l = 0$ when $2l \in q\Zset$
and $l \not \in q\Zset$,
whereas $e^-_l \neq 0$ for all $l \not \in q\Zset$.
The formulas above are equivalent to the expression of $\theta_1(t)$
given in Theorem~\ref{thm:Persistence}.

\item
We know from Lemma~\ref{lem:QR} that
$\mathcal{Q}_m = s h'_m + \tilde{\mathcal{Q}}_m$,
where $\tilde{\mathcal{Q}}_m$ is a smooth $2\pi$-periodic
function, only depending on $h_1,\ldots,h_{m-1}$
---since approximations $\varphi_{<k}$, $\psi_{<k}$ and
   $\theta_{<k}$ are uniquely determined from $h_{<k}$
   for any $k = 1, \ldots, m-1$---
that can be explicitly computed from recurrences given
in Appendix~\ref{app:Expansions}.
We also know that $\int_\Tset \mathcal{Q}_m = 0$,
and so $\int_\Tset \tilde{\mathcal{Q}}_m = 0$.
This means that there is a smooth $2\pi$-periodic function
$\zeta_m$ such that $s \zeta'_m = \tilde{\mathcal{Q}}_m$.
Finally, the necessary and sufficient condition for
$\Order(\epsilon^m)$-persistence becomes
$\mu\{ h'_m + \zeta'_m \} = 0$ or, equivalently,
$\mu_{q\Zset^\star}\{ h_m + \zeta_m \} = 0$.
\qedhere
\end{enumerate}
\end{proof}

Functions $\zeta_m$ can be recursively computed from the
formulas listed in Appendix~\ref{app:Expansions}.
For instance, we have already seen that
$\zeta_1 = 0$ and $\zeta_2 = \theta^2_1/2$ along the proof of
Theorem~\ref{thm:Persistence}.
After some tedious computations by hand
that we do not include here, we get that
\[
\zeta_3 =
\theta_1 \theta_2 + h_1 \theta_1^2/2 -
h''_1 \psi_1^2/2 + c\theta_1^3/3s - ch'_1 \theta_1 \psi_1/s,
\]
which will allow us to analyze some
$\Order(\epsilon^3)$-persistence problems in the future.
See Section~\ref{sec:OpenProblems}.
We stress that further expressions for $\zeta_4,\zeta_5,\ldots$
can be obtained using a symbolic algebra system,
but we  have doubts about their practical usefulness.

\section{Polynomial deformations of the unit circle}
\label{sec:Polynomial}

We tackle two problems in this section,
both related with polynomial deformations of the unit circle. 
First, to prove Theorem~\ref{thm:Polynomial}.
Second, to prove
that the support function~\eqref{eq:SupportFunction}
of any polynomial deformation~\eqref{eq:PolynomialPerturbation}
with $P(x,y;\epsilon) = 1 + \epsilon P_1(x,y)$ and
$P_1(x,y) \in \Rset_n[x,y]$ satisfies
condition~\eqref{eq:DegreeCondition};
see Lemma~\ref{lem:Perturbation}.
This second result allows us to explain with
Theorem~\ref{thm:Polynomial}
some of the numerical experiments about the polynomial
deformations~\eqref{eq:PolynomialPerturbation}
with $P(x,y;\epsilon) = 1 - \epsilon y^n$
performed in~\cite{MartinTamaritRamirezRos2016b},
which was the original motivation of this work.

To begin with, we check that,
once fixed the order $m \in \Nset$ and the degree $n \ge 3$,
all resonant caustics of high enough period
$\Order(\epsilon^m)$-persists under polynomial deformations
of degree $\le n$.
We explicitly quantify how high this period $q$ should be.
Some resonant caustics become more persistent in the presence
of symmetries.

\begin{theorem}
\label{thm:Degree}
Under the hypotheses of Theorem~\ref{thm:Polynomial},
the $p/q$-resonant caustic $\Order(\epsilon^m)$-persists
if any of the following three conditions are met:
\begin{enumerate}[a)]
\item
$q > nm$;
\item
$\Gamma_\epsilon$ is centrally symmetric, $q$ is odd
and $2q > nm$; or
\item
$\Gamma_\epsilon$ is anti-centrally symmetric,
$q \neq m \mod 2$, $m \ge 2$ and $q > n(m-1)$.
\end{enumerate}
\end{theorem}

\begin{lemma}
\label{lem:Perturbation}
Any deformation of the unit circle expressed
in Cartesian coordinates
as~\eqref{eq:PolynomialPerturbation} for some
$P(x,y;\epsilon) = 1 + \epsilon P_1(x,y)$ with
$P_1(x,y) \in \Rset_n[x,y]$ is a polynomial deformation
of degree $\le n$ in the sense of
Definition~\ref{defi:PolynomialDeformation}.
If $P_1$ is even (respectively, odd),
then this deformation is centrally
(respectively, anti-centrally) symmetric.
\end{lemma}

We have postponed the proofs of Theorem~\ref{thm:Degree} and
Lemma~\ref{lem:Perturbation} to Appendix~\ref{app:Proofs},
since they are technically similar.
The main difficulty is to check that some recursively computed
$2\pi$-periodic trigonometric polynomials have
degrees $\le mn$ by induction on a recursive index $m \ge 1$.

Theorem~\ref{thm:Polynomial} is a direct consequence
of Theorem~\ref{thm:Degree}.

\begin{proof}[Proof of Theorem~\ref{thm:Polynomial}]
Let $\chi = \chi(\Gamma_\epsilon,q) \in \Nset$ be
the exponent defined in Theorem~\ref{thm:Polynomial}.
Let $m = \chi - 1$.
We deduce from Theorem~\ref{thm:Degree} that the
$p/q$-resonant caustic $\Order(\epsilon^m)$-persists because:
\begin{itemize}
\item
If $\Gamma_\epsilon$ is anti-centrally symmetric,
then $nm < q$ or $n(m-1) < q$ with $m \neq q \mod 2$ and $m \ge 2$;
\item
If $\Gamma_\epsilon$ is centrally symmetric and $q$ is odd,
then $nm < 2q$; and
\item
Otherwise, $nm < q$.
\qedhere
\end{itemize}
\end{proof}

The experiments described
in~\cite[Numerical Result 5]{MartinTamaritRamirezRos2016b}
suggest that the $p/q$-resonant caustic does not
$\Order(\epsilon^\chi)$-persist under
the polynomial deformation~\eqref{eq:PolynomialPerturbation}
with $P(x,y;\epsilon) = 1 - \epsilon y^n$.
To prove it,
we should check that $h_\chi + \zeta_\chi$
has some non-zero $(q\Zset \setminus\{0\})$-harmonic.

\section{Open problems}
\label{sec:OpenProblems}

We describe three open problems that have arisen during
the development of our high-order perturbation theory.
We have not addressed them here.
Each of them is a nontrivial research challenge. 
They are work in progress.

As a general principle,
we claim that many billiard computations are greatly simplified
when working with the Bialy-Mironov generating function,
so its discovery opens the door to the resolution of
many billiard problems that seemed almost intractable.

We also stress that both our high-order perturbation theory
and our list of open problems can be extended to the setting of
dual, symplectic, wire, pensive and coin billiards,
provided we deal with those billiards in deformed circles.

\subsection{Co-preservation of resonant caustics}
\label{ssec:CoPreservation}

Tabachnikov asked if there are convex domains,
other than circles and ellipses, that possess resonant caustics
with different rotation numbers.
See~\cite[Question 4.7]{Bolsinov_etal2018}.
Theorem~\ref{thm:Persistence} and
Cyr's result on Gutkin's equation provide
an excellent framework to study the co-preservation of
resonant caustics under deformations of the unit circle. 
We mention two problems about such co-preservation.

Firstly, we consider the co-preservation of caustics with
different rotation numbers but equal periods.
Let $q \neq 2, 3, 4, 6$ be a fixed period, so
\[
\mathcal{P}_q = \{ p \in \Nset : p/q < 1/2, \ \gcd(p,q) = 1 \}
\]
has several elements.
If the deformation of the unit circle with support
function~\eqref{eq:SupportFunction} preserves
all $p/q$-resonant caustics with $p \in \mathcal{P}_q$,
then none of the smooth $2\pi$-periodic functions $h_1$ and
$h_2 + \big(\theta_1^{[p/q]}\big)^2/2$ with $p \in \mathcal{P}_q$
have resonant $(q\Zset \setminus \{ 0 \})$-harmonics, where
\[
\theta_1^{[p/q]}(t) =
\sum_{l \not \in q\Zset \cup \{-1,1\}}
\hat{\theta}_{1,l}^{[p/q]} \rme^{\rmi l t}, \qquad
\hat{\theta}_{1,l}^{[p/q]} = \nu_l(p/q) \hat{h}_{1,l},
\]
for some functions $\nu_l:(0,1/2) \to \Rset$ defined
in~\eqref{eq:Factors},
which have no rational roots when $|l| \ge 2$.
The above necessary conditions for $\Order(\epsilon^2)$-persistence
impose rather stringent restrictions on
the first-order coefficient $h_1$.
We do not state any specific result because we have not
yet found one satisfactory enough.

Secondly, we consider the co-preservation of caustics
with rotation numbers $1/2$ and $1/q$ for some fixed period
$q \ge 3$ under deformations of the unit circle.
The preservation of the $1/2$-resonant caustic implies that
all $(2\Zset \setminus\{0\})$-harmonics of the
support function~\eqref{eq:SupportFunction} are zero,
which greatly simplifies the problem.
The case $q=3$ for centrally symmetric deformations was studied by
J.~Zhang~\cite{Zhang2019}.
His method consists in a careful analysis of
the obstructions for $\Order(\epsilon^2)$-persistence of
the $1/3$-resonant caustic.
His computations are more involved
than the ones outlined here because he uses the standard
generating function for billiards (that is,
the minus distance between consecutive impacts) and
he considers deformations of the unit circle written
in polar coordinates.
We plan to extend Zhang's results to some periods $q>3$,
starting with periods $q=4,5$ for which it is necessary to
study the obstructions for $\Order(\epsilon^3)$-persistence.

\subsection{A convergence problem}
\label{ssec:Convergence}

We recall that the function $\zeta_m$ in Theorem~\ref{thm:Persistence} only depends on $h_1,\ldots,h_{m-1}$.
Hence,
if the $p/q$-resonant caustic $\Order(\epsilon^m)$-persists,
but not $\Order(\epsilon^{m+1})$-persists,
and we fix any higher order $r > m$,
then there is a smooth $2\pi$-periodic function
$\eta^{[r]}_\epsilon = \Order(\epsilon^{m+1})$
such that the resonant caustic $\Order(\epsilon^r)$-persists
under the new deformation
$\Gamma^{[r]}_\epsilon =
 \Gamma_\epsilon + \Order(\epsilon^{m+1})$
with support function
$h^{[r]}_\epsilon  = h_\epsilon + \eta^{[r]}_\epsilon$.
This function $\eta^{[r]}_\epsilon$ can also be
recursively computed and
only contains $(q\Zset \setminus \{0\})$-harmonics.
See the last item in Proposition~\ref{prop:Iterative}.

A natural question is:
Does $h^{[\infty]}_\epsilon :=
 \lim_{r \to +\infty} h^{[r]}_\epsilon$ exist?
That is, we look for a corrected support function
$h^{[\infty]}_\epsilon = h_\epsilon + \Order(\epsilon^{m+1})$
such that the $p/q$-resonant caustic persists \emph{at all orders}
under the deformation $\Gamma^{[\infty]}_\epsilon$
with support function $h^{[\infty]}_\epsilon$.
This problem is closely related to the density property
established in~\cite{KaloshinZhang2018}.

The first author, in collaboration with V.~Kaloshin and K.~Zhang,
is working on a more ambitious version of this problem.
Namely,
once \emph{fixed} a rational rotation number $1/q$,
the goal is to construct a functional
\[
\mathcal{H}_{\Zset \setminus q\Zset} \ni h
\stackrel{\mathcal{F}_{1/q}}{\mapsto}
\eta \in \mathcal{H}_{q\Zset \setminus\{0\}}
\]
such that:
1) $\mathcal{H}_A$ is a neighborhood of zero in a suitable space
of $2\pi$-periodic functions with only $A$-harmonics;
2) $\mathcal{F}_{1/q}(0) = 0$;
3) $\mathcal{F}_{1/q}$ is as regular as possible at $h = 0$,
and 4) the convex domain with support function
$1 + h + \mathcal{F}_{1/q}[h]$ has a $1/q$-resonant caustic
for any $h \in \mathcal{H}_{\Zset \setminus q\Zset}$.

\subsection{Exponentially small phenomena}
\label{ssec:ExponentiallySmall}

Let $\Gamma_\epsilon$ be a polynomial deformation
of the unit circle of degree $\le n$ in the sense of
Definition~\ref{defi:PolynomialDeformation}.
Let $q \ge 3$ be a \emph{fixed} period.
Theorem~\ref{thm:Polynomial} implies that
the $p/q$-resonant caustics $\Order(\epsilon^{\chi - 1})$-persist under $\Gamma_\epsilon$ for all $p \in \mathcal{P}_q$,
so there are perturbed normal functions $\psi^{[p/q]}_\epsilon(t)$
and perturbed incidence-reflection functions
$\theta^{[p/q]}_\epsilon(t)$ such that the error
\[
\mathcal{E}^{[p/q]}_\epsilon(t) :=
\sigma
\big\{
h'_\epsilon(\psi^{[p/q]}_\epsilon(t)) \cdot
\sin \theta^{[p/q]}_\epsilon(t)
\big\} -
\delta
\big\{
h_\epsilon(\psi^{[p/q]}_\epsilon(t)) \cdot
\cos \theta^{[p/q]}_\epsilon(t)
\big\}
\]
is $\Order(\epsilon^\chi)$ as $\epsilon \to 0$.
We recall that $\chi \asymp 2q/n$ as odd $q \to +\infty$
for centrally symmetric deformations,
and $\chi \asymp q/n$ as $q \to +\infty$ otherwise.
Let us focus on the second case, which is the generic one.
It is natural to ask whether, in this generic case,
\[
\mathcal{E}^{[p/q]}_\epsilon(t) =
\Order( \epsilon^\chi ) \simeq
\Order( \epsilon^{q/n} ) =
\Order( \rme^{-q |\log \epsilon|/n} )
\]
as $q \to +\infty$ for \emph{fixed}
$\epsilon \in [-\epsilon_0,\epsilon_0]$.
This is a hard problem.
It was partially addressed in~\cite{MartinTamaritRamirezRos2016a},
where the authors establish an exponentially small upper bound
on the difference of lengths of $1/q$-periodic billiard
trajectories in analytic strictly convex domains as $q \to +\infty$.
The numerical experiments described
in~\cite{MartinTamaritRamirezRos2016b} suggest that these
upper bounds can be improved to exponentially small
asymptotic formulas for some deformations.
To be precise,
we deduce from the computations about the polynomial
deformations~\eqref{eq:PolynomialPerturbation} with
$P(x,y;\epsilon) = 1 - \epsilon y^n$,
that there are constants $\xi(\Gamma_\epsilon)$
with a finite limit as $\epsilon \to 0$ such that the error
$\mathcal{E}^{[1/q]}_\epsilon(t)$ has `size'
$q^{-3} \rme^{-q [|\log \epsilon|/n) + \xi(\Gamma_\epsilon)]}$
as $q \to +\infty$ for any fixed
$\epsilon \in [-\epsilon_0,\epsilon_0]$.

This problem can be addressed by
the direct approach used by Wang in~\cite{Wang2020} or
the approach used by Kaloshin and Zhang
in~\cite{KaloshinZhang2018}.
However,
we feel that exponentially small asymptotic formulas
can only be obtained with more refined techniques,
like the extension of $\psi^{[p/q]}_\epsilon(t)$ and
$\theta^{[p/q]}_\epsilon(t)$ to complex values of $t$,
the analysis of their complex singularities,
resurgence theory and complex matching.
See~\cite{DelshamsRamirez1998,Gelfreich1999,MartinSauzinSeara2011}
for examples of how these refined techniques are applied
in the setting of discrete systems (maps).

\section*{Acknowledgments}
C.~E.~K. gratefully acknowledges support from
the European Research Council (ERC) through
the Advanced Grant “SPERIG” (\#885 707) .
R.~R.-R. was supported in part by the grant PID-2021-122954NB-100
which was funded by MCIN/AEI/10.13039/501100011033 and
“ERDF: A way of making Europe”.
R.~R.-R. thanks to  Pau Mart{\'\i}n and Tere Seara for
many long and stimulating conversations on related problems.
We are also indebted to the referee for 
suggestions that helped us to improve the exposition.

\appendix

\section{Shift, sum, difference and average operators}
\label{app:Operators}

Let $0 < p < q$ integers such that $\gcd(p,q) = 1$.
Set $\omega = 2\pi p/q$.
The shift, sum, difference and average operators act on smooth
functions $a: \Tset \to \Rset$ as follows
\begin{align*}
\tau\{ a(t) \} &= a(t + \omega),\\
\sigma\{ a(t) \} &= a(t + \omega) + a(t), \\
\delta\{ a(t) \} &= a(t + \omega) - a(t), \\
\mu \{ a(t) \} &= \textstyle \frac{1}{q} \sum_{j=0}^{q-1} a(t + j \omega).
\end{align*}
They diagonalize in the Fourier basis and the average operator
is the projection onto the resonant $q\Zset$-harmonics.

\begin{lemma}
\label{lem:Fourier}
If $a(t) = \sum_{l\in\Zset} \hat{a}_l \rme^{\rmi l t}$
is a smooth $2\pi$-periodic function, then
\begin{align*}
\tau\{ a(t) \} &=
\textstyle \sum_{l \in \Zset} \rme^{\rmi l \omega} \hat{a}_l \rme^{\rmi l t},\\
\sigma\{ a(t) \} &=
\textstyle \sum_{l \in \Zset} (\rme^{\rmi l \omega} + 1 )\hat{a}_l \rme^{\rmi l t}, \\
\delta\{ a(t) \} &=
\textstyle \sum_{l \in \Zset} (\rme^{\rmi l \omega} - 1) \hat{a}_l \rme^{\rmi l t}, \\
\mu \{ a(t) \} &=
\textstyle \sum_{l \in q\Zset} \hat{a}_l \rme^{\rmi l t}.
\end{align*}
Hence,
$\mu \circ \tau = \mu \circ \mu = \mu$,
$\mu \circ \sigma = 2\mu$ and $\mu \circ \delta = 0$
on the set of smooth $2\pi$-periodic functions.
\end{lemma}

\begin{proof}
The average $\frac{1}{q}\sum_{j=0}^{q-1} \rme^{\rmi j l \omega}$ is equal
to one when $l \in q\Zset$ but it vanishes otherwise.
\end{proof}

We deduce from this lemma that
persistence condition~\eqref{eq:PersistenceCondition} holds
if and only if all resonant $q\Zset$-harmonics of
$\mathcal{Q}_k$ are equal to zero.

Next, we invert operator $\delta$,
which is the key point in solving~\eqref{eq:CompatibleEquations}
and~\eqref{eq:RelationsCoefficients}.

\begin{lemma}
\label{lem:Inversion}
Let $b(t)$ be any smooth $2\pi$-periodic function.
Equation
\[
\delta\{ a(t) \} = b(t)
\]
has a smooth $2\pi$-periodic solution $a(t)$
if and only if $\mu\{ b(t) \} = 0$,
in which case the solution is unique once we fix
its smooth $2\pi/q$-periodic average $\mu\{a(t)\}$.
Analogously, equation
\[
\delta\{ t + \tilde{a}(t) \} = b(t)
\]
has a smooth $2\pi$-periodic solution $\tilde{a}(t)$
if and only if $\mu\{ b(t) \} \equiv \omega$,
in which case the solution is unique once we fix
its smooth $2\pi/q$-periodic average $\mu\{\tilde{a}(t)\}$.
\end{lemma}

\begin{proof}
If $a(t)$ is a $2\pi$-periodic solution, then
\[
\mu \{ b(t) \} =
\mu \big \{ a(t + \omega) - a(t) \big \} =
\frac{1}{q} \sum_{j=0}^{q-1} [a(t+(j+1)\omega) - a(t + j\omega) ] =
\frac{a(t + q\omega) - a(t) }{q} = 0.
\]
Conversely, if $\mu\{ b(t) \} = 0$,
then $b(t) = \sum_{l \not \in q\Zset} \hat{b}_l \rme^{\rmi l t}$
so that
\[
\hat{a}_l = \hat{b}_l/(\rme^{\rmi l \omega} - 1), \qquad
\forall l \not \in q\Zset.
\]
defines the Fourier coefficients of the unique
smooth $2\pi$-periodic solution
$a(t) = \sum_{l \not \in q\Zset} \hat{a}_l \rme^{\rmi l t}$
such that $\mu \{ a(t) \} = 0$.
Smoothness follows from the bounds $|\hat{a}|_l \le \nu |\hat{b}_l|$
for all $l \not \in q \Zset$, where $\nu = 1/|\rme^{\rmi 2 \pi/q} -1| > 0$.
Since $\delta \circ \mu = 0$,
we deduce that the average $\mu\{ a(t) \}$ does not need to be zero,
but can be determined freely.

If $a(t) = t + \tilde{a}(t)$ is a solution
such that $\tilde{a}(t)$ is $2\pi$-periodic, then
\[
\mu \{ b(t) \} =
\frac{a(t + q\omega) - a(t) }{q} =
\frac{ q\omega + \tilde{a}(t + q\omega) - \tilde{a}(t)}{q} = \omega.
\]
The converse is obtained with the same argument as before.
\end{proof}

\section{Computation of some asymptotic expansions}
\label{app:Expansions}

Recall that $\omega = 2\pi p/q$, $c = \cos(\omega/2)$ and
$s = \sin(\omega/2) > 0$, where $p/q \in (0,1/2)$ is
some rational number such that $\gcd(p,q) = 1$.
Let us compute the asymptotic expansions of
\[
\mathcal{R} := h \circ \psi \cdot \cos \theta, \qquad
\mathcal{Q} := h' \circ \psi \cdot \sin \theta,
\]
from the asymptotic expansion~\eqref{eq:SupportFunction}
of the support function and the
asymptotic expansions~\eqref{eq:FunctionsExpasions}
of the side, normal and incidence-reflection functions.
We assume that all coefficients $h_k(\psi)$, $\varphi_k(t)$,
$\psi_k(t)$ and $\theta_k(t)$ are smooth $2\pi$-periodic functions.

\begin{lemma}
\label{lem:CS}
The coefficients of the asymptotic expansions
\[
\sin \theta \asymp \sum_{k \ge 0} \mathcal{S}_k \epsilon^k,\quad
\cos \theta \asymp \sum_{k \ge 0} \mathcal{C}_k \epsilon^k
\quad \mbox{ as } \epsilon \to 0
\]
can be computed from the initial values
$\mathcal{S}_0 = s$ and $\mathcal{C}_0 = c$
by means of the recurrences
\begin{equation}
\label{eq:CS}
\mathcal{S}_k =
c \theta_k + \frac{1}{k} \sum_{l=1}^{k-1} l \theta_l \mathcal{C}_{k-l},
\quad
\mathcal{C}_k =
-s \theta_k - \frac{1}{k} \sum_{l=1}^{k-1} l \theta_l \mathcal{S}_{k-l},
\quad
\forall k \ge 1.
\end{equation}
In particular,
$\mathcal{S}_k =
 \mathcal{S}_k \llbracket \theta_{\le k} \rrbracket =
 c \theta_k + \tilde{\mathcal{S}}_k \llbracket \theta_{< k} \rrbracket$ and
$\mathcal{C}_k =
 \mathcal{C}_k \llbracket \theta_{\le k} \rrbracket =
 -s \theta_k + \tilde{\mathcal{C}}_k \llbracket \theta_{< k} \rrbracket$.
\end{lemma}

\begin{proof}
These recurrences are directly obtained from identities
\[
\frac{\rmd}{\rmd \epsilon} \{ \sin \theta \} =
\frac{\rmd \theta}{\rmd \epsilon} \cdot \cos \theta, \quad
\frac{\rmd}{\rmd \epsilon} \{ \cos \theta \} =
-\frac{\rmd \theta}{\rmd \epsilon} \cdot \sin \theta.
\qedhere
\]
\end{proof}

\begin{definition}
If $\alphav = \{\alpha_1,\alpha_2,\ldots\}$ is a sequence of
non-negative integers with a finite number of non-zero terms
and $\psiv = \{\psi_1,\psi_2,\ldots\}$ is a sequence of
smooth $2\pi$-periodic functions, then
\[
|\alphav| = \sum_{l \ge 1} \alpha_l, \quad
\| \alphav \| = \sum_{l \ge 1} l\alpha_l, \quad
\alphav! = \prod_{l \ge 1} \alpha_l!, \quad
\psiv^\alphav = \prod_{l \ge 1} \psi_l^{\alpha_l}.
\]
\end{definition}

\begin{lemma}
\label{lem:H}
The coefficients of the asymptotic expansion
\[
h \circ \psi \asymp 1 + \sum_{k \ge 1} \mathcal{H}_k \epsilon^k
\quad \mbox{ as } \epsilon \to 0
\]
are given by
$\mathcal{H}_k =
 \mathcal{H}_k \llbracket h_{\le k},\psi_{<k} \rrbracket =
 h_k + \tilde{\mathcal{H}}_k \llbracket h_{<k},\psi_{<k} \rrbracket$ with
\begin{equation}
\label{eq:H}
\tilde{\mathcal{H}}_k =
\sum_{i=1}^{k-1} \sum_{j=1}^{k-i}
\left(\sum_{|\alphav| = j,\| \alphav \| = k-i}
\frac{\psiv^\alphav}{\alphav !} \right) h_i^{(j)}, \quad
\forall k \ge 1.
\end{equation}
Analogously,
$h' \circ \psi \asymp
 \sum_{k \ge 1}
 \mathcal{H}_k\llbracket h'_{\le k},\psi_{<k} \rrbracket \epsilon^k$
as $\epsilon \to 0$.
\end{lemma}

\begin{proof}
Set $\Delta \psi := \psi - t \asymp \sum_{i \ge 1} \epsilon^i \psi_i$
as $\epsilon \to 0$.
Formula~\eqref{eq:H} is a direct consequence of the asymptotic expansions
\[
h_i(t + \Delta \psi) \asymp
\sum_{j \ge 0} \frac{(\Delta \psi)^j}{j!} h_i^{(j)},\qquad
\frac{(\Delta \psi)^j}{j!} \asymp
\sum_{|\alphav| = j} \frac{\psiv^\alphav}{\alphav!} \epsilon^{\| \alphav \|}
\]
as $\Delta \psi \to 0$ and $\epsilon \to 0$, respectively.
The first expansion is the Taylor theorem.
The second expansion is the multinomial theorem.

On the one hand,
if $| \alphav | = j = 0$ and $\| \alphav \| = k-i$,
then $i=k$ and $\alphav = \{0,0,\ldots\}$.
On the other hand,
if $| \alphav | = j \ge 1$ and $\| \alphav \| = k - i$,
then $i \le k-1$, $j \le k - i$ and $\alpha_l = 0$ for all $l \ge k$.
This justifies that
$\mathcal{H}_k =
 \mathcal{H}_k \llbracket h_{\le k},\psi_{<k} \rrbracket =
 h_k + \tilde{\mathcal{H}}_k \llbracket h_{<k},\psi_{<k} \rrbracket$.
\end{proof}

The next result is a more informative version of Lemma~\ref{lem:QR},
whose proof was pending.

\begin{lemma}
\label{lem:QRbis}
The coefficients of the asymptotic expansions
\[
\mathcal{R} = h \circ \psi \cdot \cos \theta \asymp
c + \sum_{k \ge 1} \mathcal{R}_k \epsilon^k, \quad
\mathcal{Q} = h' \circ \psi \cdot \sin \theta \asymp
\sum_{k \ge 1} \mathcal{Q}_k \epsilon^k
\quad \mbox{ as } \epsilon \to 0
\]
are given by
\begin{align*}
\mathcal{Q}_k &=
\mathcal{Q}_k \llbracket h'_{\le k},\psi_{<k},\theta_{<k} \rrbracket =
s h'_k +
\tilde{\mathcal{Q}}_k \llbracket h'_{<k},\psi_{<k},\theta_{<k} \rrbracket, \\
\tilde{\mathcal{Q}}_k &=
s \tilde{\mathcal{H}}_k \llbracket (h'_{<k},\psi_{<k} \rrbracket +
\textstyle \sum_{l=1}^{k-1}
\mathcal{H}_l \llbracket h'_{\le l},\psi_{<l} \rrbracket
\mathcal{S}_{k-l} \llbracket \theta_{\le k-l} \rrbracket, \\
\mathcal{R}_k &=
\mathcal{R}_k \llbracket h_{\le k},\psi_{<k},\theta_{\le k} \rrbracket =
\tilde{\mathcal{R}}_k \llbracket h_{\le k},\psi_{<k},\theta_{<k} \rrbracket -
 s\theta_k, \\
\tilde{\mathcal{R}}_k &=
c \mathcal{H}_k \llbracket h_{\le k},\psi_{<k} \rrbracket +
\tilde{\mathcal{C}}_k \llbracket \theta_{<k} \rrbracket +
\textstyle \sum_{l=1}^{k-1}
\mathcal{H}_l \llbracket h_{\le l},\psi_{<l} \rrbracket
\mathcal{C}_{k-l} \llbracket \theta_{\le k-l} \rrbracket.
\end{align*}
\end{lemma}

\begin{proof}
It is a direct consequence of Lemmas~\ref{lem:CS} and~\ref{lem:H}.
\end{proof}

Finally, we compute the asymptotic expansion of the support function
of the deformation of the unit circle given
by~\eqref{eq:PolynomialPerturbation} in Cartesian coordinates.
We assume that $P(x,y;\epsilon) = 1 + \epsilon P_1(x,y)$ and
$P_1(x,y) \in \Rset_n[x,y]$.

\begin{lemma}
\label{lem:G}
If $P(x,y;\epsilon) = 1 + \epsilon P_1(x,y)$ with
$P_1(x,y) = \sum_{i,j \ge 0, i+j \le n} p_{ij} x^i y^j$,
then the coefficients of the asymptotic
expansion~\eqref{eq:SupportFunction} of the support function of
the deformation~\eqref{eq:PolynomialPerturbation}
can be computed from recurrences
\begin{equation}
\label{eq:hRecurrences}
2h_k + \tilde{\mathcal{G}}^\star_k \llbracket h_{<k} \rrbracket +
\tilde{\mathcal{G}}^\bullet_k \llbracket h'_{<k} \rrbracket =
\mathcal{G}^\diamond_{k-1}
\llbracket h_{\le k-1},h'_{\le k-1} \rrbracket,
\qquad \forall k \ge 1,
\end{equation}
where
\begin{align*}
\mathcal{G}^\star_k &=
2h_k + \tilde{\mathcal{G}}^\star_k =
2h_k + \textstyle 2 \sum'_{|\alphav| = 2, \| \alphav \| = k} \hv^\alphav/\alphav!, \\
\tilde{\mathcal{G}}^\bullet_k &=
\textstyle 2 \sum'_{|\alphav| = 2, \| \alphav \| = k} (\hv')^\alphav/\alphav!, \\
\mathcal{G}^\diamond_k &=
\textstyle \sum_{i,j \ge 0, i+j \le n} p_{ij} \ i! \ j!
\sum_{|\alphav| = i, |\betav| = j, \| \alphav \| + \| \betav \| = k}
\xv^\alphav \yv^\betav/ \alphav! \betav!, \\
\xv &= \cos \psi \cdot \hv - \sin \psi \cdot \hv', \\
\yv &= \sin \psi \cdot \hv + \cos \psi \cdot \hv'.
\end{align*}
Here, $\alphav = \{ \alpha_0, \alpha_1, \ldots\}$,
$\betav = \{ \beta_0, \beta_1, \ldots\}$,
$\hv = \{1, h_1, h_2, \ldots\}$ and
$\hv' = \{0, h'_1, h'_2,\ldots\}$.
Symbol $\sum'_{|\alphav| = 2, \| \alphav \| = k}$ means that we do not
include the term with $\alpha_0 = 1$ and $\alpha_k = 1$.
\end{lemma}

\begin{proof}
We know that
$z(\psi;\epsilon) = \big( x(\psi;\epsilon), y(\psi;\epsilon) \big)$,
where
\begin{align*}
x = x(\psi;\epsilon) &=
\cos \psi \cdot h(\psi;\epsilon) - \sin \psi \cdot h'(\psi;\epsilon),\\
y = y(\psi;\epsilon) &=
\sin \psi \cdot h(\psi;\epsilon) + \cos \psi \cdot h'(\psi;\epsilon),
\end{align*}
is a normal parametrization of $\Gamma_\epsilon$,
so the support function satisfies the implicit equation
\[
h^2 + (h')^2 = x^2 + y^2 = 1 + \epsilon P_1(x,y).
\]
Therefore, recurrence~\eqref{eq:hRecurrences} is a direct consequence of
the asymptotic expansions
\begin{align*}
h^2 \asymp
1 + \sum_{k \ge 1} \mathcal{G}^\star_k \epsilon^k,\quad
(h')^2 \asymp
\sum_{k \ge 2} \tilde{\mathcal{G}}^\bullet_k \epsilon^k, \quad
P_1(x,y) \asymp
\sum_{k \ge 0} \mathcal{G}^\diamond_k \epsilon^k
\end{align*}
as $\epsilon \to 0$,
all of which follow from the multinomial theorem.
\end{proof}

\section{Proofs of Theorem~\ref{thm:Degree} and
Lemma~\ref{lem:Perturbation}}
\label{app:Proofs}

There are two main tools for both proofs.
Firstly, the explicit formulas for the asymptotic coefficients
$\mathcal{S}_k = c \theta_k + \tilde{\mathcal{S}}_k$,
$\mathcal{C}_k = \tilde{\mathcal{C}}_k - s\theta_k$,
$\mathcal{H}_k = h_k + \tilde{\mathcal{H}}_k$,
$\mathcal{R}_k = \tilde{\mathcal{R}}_k - s\theta_k$,
$\mathcal{Q}_k$,
$\tilde{\mathcal{G}}^\star_k$,
$\tilde{\mathcal{G}}^\bullet_k$ and
$\mathcal{G}^\diamond_k$
given in Appendix~\ref{app:Expansions}.
Secondly, the following elementary properties:
\begin{itemize}
\item
$T_k[t]$ is a real vector space;
\item
$a(t) \in T_k[t], \ b(t) \in T_l[t] \Rightarrow
a(t)b(t) \in T_{k+l}[t], \ a'(t) \in T_k[t]$;
\item
$\alphav = \{ \alpha_0,\alpha_1,\ldots \}$ and
$\av = \{ a_0,a_1,\ldots\}$ with
$a_j(t) \in T_{nj}[t] \Rightarrow
\av^\alphav(t) \in T_{n \| \alphav \|}[t]$;
\item
$a(t) \in T_k[t] \Rightarrow
\sigma\{a(t)\}, \delta\{ a(t) \}, \mu\{ a(t) \} \in T_k[t]$; and
\item
$b(t) \in T_k[t], \ \mu\{ b(t) \} = 0 \Rightarrow
\exists ! a(t) \in T_k[t]$ s. t. $\delta\{ a(t) \} = b(t)$
and $\mu\{ a(t) \} = 0$.
\end{itemize}
We will use these properties without any explicit mention in what follows.

\begin{lemma}
\label{lem:GeneralCase}
If condition~\eqref{eq:DegreeCondition} holds and $q > nm$,
then we can compute the $\Order(\epsilon^m)$-corrections
$\theta_{\le m}(t;\epsilon)$, $\varphi_{\le m}(t;\epsilon)$ and
$\psi_{\le m}(t;\epsilon)$ by solving the compatible
equations~\eqref{eq:CompatibleEquations}
and~\eqref{eq:RelationsCoefficients} for $k=1,\ldots,m$.
Compatibility is guaranteed because all necessary and sufficient
persistence conditions~\eqref{eq:PersistenceCondition} hold
for $k=1,\ldots,m$.
\end{lemma}

\begin{proof}
It suffices to check that
\begin{enumerate}[i)]
\item
$\mathcal{H}_k = h_k + \tilde{\mathcal{H}}_k \in T_{nk}[t]$;
\item
$\tilde{\mathcal{R}}_k,\mathcal{Q}_k \in T_{nk}[t]$
---so condition~\eqref{eq:PersistenceCondition} holds because
   $q > nm$ and $\bar{\mathcal{Q}}_k = 0$---;
\item
$\theta_k, \varphi_k, \psi_k \in T_{nk}[t]$; and
\item
$\mathcal{S}_k, \mathcal{C}_k \in T_{nk}[t]$;
\end{enumerate}
for $k=1,\ldots,m$.
We prove it by induction over $m$.

The base case $m=1$ is trivial.
Namely,
$\mathcal{H}_1 = h_1 \in T_n[t]$ by
condition~\eqref{eq:DegreeCondition},
$\tilde{\mathcal{R}}_1 = ch_1 \in T_n[t]$,
$\mathcal{Q}_1 = sh'_1 \in T_n[t]$,
$\theta_1 \in T_n[t]$ is the unique solution of
problem~\eqref{eq:CompatibleEquation1},
$\varphi_1 \in T_n[t]$ is the unique solution of
equation $\delta\{ \varphi_1 \} = 2\theta_1$ such that $\mu\{ \varphi_1 \} = 0$,
$\psi_1 = \varphi_1 + \theta_1 \in T_n[t]$,
$\mathcal{S}_1 = c \theta_1 \in T_n[t]$ and
$\mathcal{C}_1 = -s \theta_1 \in T_n[t]$.

Next, let us assume that properties i)--iv) hold for
$k=1,\ldots,m-1$.
We need to prove that they also hold for $k=m$.
Property~i) follows from~\eqref{eq:H} and
condition~\eqref{eq:DegreeCondition}.
Property~ii) follows from the recurrences given in Lemma~\ref{lem:QRbis}.
Then $\theta_m \in T_{nm}[t]$ is the unique solution of equation
$s \delta\{ \theta_m \} =
 \delta\{ \tilde{\mathcal{R}}_m \} -\sigma\{ \mathcal{Q}_m \}$
such that $\mu\{ \theta_m \} = 0$,
$\varphi_m \in T_{nm}[t]$ is the unique solution of
equation $\delta\{ \varphi_m \} = 2\theta_m$
such that $\mu\{ \varphi_m \} = 0$,
and $\psi_m = \varphi_m + \theta_m \in T_{nm}[t]$.
This proves property~iii).
Property~iv) follows from recurrences~\eqref{eq:CS}.
\qedhere
\end{proof}

\begin{lemma}
If condition~\eqref{eq:DegreeCondition} holds,
$\Gamma_\epsilon$ is centrally symmetric, $q$ is odd and $2q > nm$,
then we can compute all $\Order(\epsilon^m)$-corrections too.
\end{lemma}

\begin{proof}
If $\Gamma_\epsilon$ is centrally symmetric,
then its support function $h(\psi;\epsilon)$ is $\pi$-periodic in $\psi$.
In that case, it is not hard to prove by induction over $m$ that
objects~i)--iv) listed in the proof of the previous lemma
are also $\pi$-periodic for $k=1,\ldots,m$.
In particular, if $q$ is odd and $2q > nm$,
then persistence condition~\eqref{eq:PersistenceCondition} holds
because all resonant $q\Zset$-harmonics of the
$\pi$-periodic trigonometric polynomial $\mathcal{Q}_k \in T_{nk}[t]$
are equal to zero for $k=1,\ldots,m$.
\end{proof}

\begin{lemma}
If condition~\eqref{eq:DegreeCondition} holds,
$\Gamma_\epsilon$ is anti-centrally symmetric,
$q \neq m \mod 2$, $m \ge 2$ and $q > n(m-1)$,
then we can compute all $\Order(\epsilon^m)$-corrections too.
\end{lemma}

\begin{proof}
We already know from Lemma~\ref{lem:GeneralCase} that
we can compute the $\Order(\epsilon^{m-1})$-corrections,
since $q > n(m-1)$.
Therefore, we only need to check that the last persistence condition
\[
\mu\{ \mathcal{Q}_m \} = 0
\]
holds.
That is, we need to check that all resonant
$q \Zset$-harmonics of $\mathcal{Q}_m$ are equal to zero.

If the perturbation $\Gamma_\epsilon$ is anti-centrally symmetric,
then its support function satisfies that
$h(\psi;\epsilon) = h(\psi+\pi;-\epsilon)$ and
its asymptotic coefficients satisfy that
\[
h_k(\psi + \pi) = (-1)^k h_k(\psi), \quad \forall k \in \Nset.
\]
In that case, it is not hard to prove by induction over $m$
that objects~i)--iv) listed in the proof of
Lemma~\ref{lem:GeneralCase} satisfy the same property.
Namely, that they are $\pi$-periodic and $\pi$-antiperiodic functions
for even and odd indexes $k$, respectively.
In particular,
if $m$ is even (respectively, odd),
then $\mathcal{Q}_m$ only contains even (respectively, odd) harmonics.
Hence,
$\mathcal{Q}_m$ does not contain harmonics
$\rme^{\pm \rmi q t}$, because $q \neq m \mod 2$,
and it does not contain any harmonic
$\rme^{\pm \rmi l q t}$ with $l \ge 2$ either,
because $lq \ge 2q > 2n(m-1) \ge nm$ and
$\mathcal{Q}_m \in T_{nm}[t]$.
We have used that $m \ge 2$ in the last inequality.
\end{proof}

Theorem~\ref{thm:Degree} is a direct consequence
of the previous three lemmas.

\begin{proof}[Proof of Lemma~\ref{lem:Perturbation}]
We have to check that
$\tilde{\mathcal{G}}^\star_k,
\tilde{\mathcal{G}}^\bullet_k,
\mathcal{G}^\diamond_{k-1}, h_k \in T_{nk}$
for all $k \ge 1$.
We prove it by induction over $k$.
The base case $k=1$ is trivial,
because $\tilde{\mathcal{G}}^\star_1 = 0$,
$\tilde{\mathcal{G}}^\bullet_1 = 0$,
$\mathcal{G}^\diamond_0 = P_1(\cos \psi,\sin \psi)$
and $h_1 = \frac{1}{2} P_1(\cos \psi,\sin \psi)$.
The induction step follows from the explicit formulas given in
Lemma~\ref{lem:G} and the elementary properties listed
at the beginning of this appendix.
The claims about symmetries are trivial.
\end{proof}

\end{document}